\documentclass[11pt, reqno]{amsart}
\usepackage[margin=2.75cm]{geometry}

\usepackage{amsthm,
amsmath,mathrsfs,enumitem,
amssymb,mathtools,mathrsfs,
dsfont,centernot,graphicx,
xcolor,hyperref,tikz,subcaption}
\usetikzlibrary{decorations.pathreplacing,calligraphy}

\usepackage[english]{babel}
\usepackage[hypcap]{caption}

\theoremstyle{plain}
\newtheorem{theorem}{Theorem}[section]

\theoremstyle{definition}
\newtheorem{definition}[theorem]{Definition}

\newtheorem{assumption}[theorem]{Assumption}
\newtheorem{remark}[theorem]{Remark}
\numberwithin{equation}{section}

\newcommand{\indep}{\perp \!\!\! \perp}
\newcommand{\Leb}{\mathrm{Leb}}

\newcommand{\rrvert}{\vert}
\newcommand{\llvert}{\vert}

\numberwithin{equation}{section}
\makeatother

\usepackage{graphicx}
\usepackage{mathrsfs}
\usepackage{tikz}
\usetikzlibrary{trees,snakes,external,backgrounds}

\title[Interacting processes on sparse random graphs]{Interacting stochastic processes on sparse random graphs}\thanks{The author was  supported by the Office of Naval Research under the Vannevar Bush Faculty Fellowship
  N0014-21-1-2887, and  the National Science Foundation under grant
  DMS-1954351} 
\author{Kavita Ramanan}
\address{Division of Applied Mathematics, Brown University, Providence, RI 02912, USA}
\email{kavita\_ramanan@brown.edu}
\begin{document}
\maketitle

\begin{abstract}
Large ensembles of stochastically evolving interacting particles
describe phenomena in diverse fields including statistical physics, neuroscience, biology, and engineering.
In such systems, the infinitesimal evolution of each particle depends only on its own state (or history)
and the states (or histories) of neighboring particles with respect to
an underlying, possibly random, interaction graph.  While these high-dimensional processes are typically too
complex to be amenable to exact analysis, their dynamics are quite well understood when the interaction graph is
the complete graph. In this case, classical theorems show that in the limit as the number of particles goes to infinity,  the dynamics of the empirical measure and the law of a typical particle coincide and can be characterized in terms of a much more tractable dynamical system of reduced dimension called the mean-field limit.  In contrast,  until recently not much was known about corresponding convergence results in the complementary case when the interaction graph is sparse (i.e., with uniformly bounded average degree). This article provides a brief survey of classical
work and then describes recent progress on the sparse regime that relies
on a combination of techniques from random graph theory, Markov random
fields, and stochastic analysis. The article concludes by discussing ramifications
for applications and posing several open problems.
\end{abstract}



\maketitle

\section{Introduction}
\label{188.sec1}

\subsection{Background}
\label{subs-back}

A recurring theme in probability theory is the emergence of deterministic
(or more predictable) behavior when there is an aggregation of many random
elements. A classical result is the strong law of large numbers established
by Kolmogorov in 1933 \cite{Kol-SLLN33}. This  states that given a sequence
of random variables $(X_{i})_{i \in {\mathbb{N}}}$ that are independent
and identically distributed (i.i.d.) and have finite mean (equivalently,
$(X_{i})_{i \in {\mathbb{N}}}$ is distributed according to some product
probability measure ${\otimes }^{{\mathbb{N}}}\nu $, where $\nu $ is a probability
measure on the Borel sets of $\mathbb{R}$ that satisfies
$\int_{\mathbb{R}} |x|  \nu (dx) < \infty $), then with probability
one,

\begin{equation}
\label{SLLN} S_{n} := \frac{1}{n} \sum
_{i=1}^{n} X_{i} \rightarrow {
\mathbb{E}}[X_{1}] = \int_{\mathbb{R}} x \nu (dx), \quad
\mbox{as } n \rightarrow \infty .
\end{equation}
In a similar spirit, the Glivenko--Cantelli theorem, also established in
1933 \cite{Gli-33,Cantelli-33b}, provides information on the asymptotic behavior of empirical
measures of i.i.d. random variables.  Specifically, it  shows that
 with probability one,

\begin{equation}
\label{G-C} \mu _{n} := \frac{1}{n} \sum
_{i=1}^{n} \delta _{X_{i}} \rightarrow {
\mathbb{E}}[ \delta _{X_{1}}] = \mathscr{L}(X_{1}) = \nu ,
\quad \mbox{as } n \rightarrow \infty ,
\end{equation}
where $\delta _{x}$ represents the Dirac delta measure at $x$ and
$\mathscr{L}(Y)$ denotes the law or distribution of a random variable
$Y$.  The convergence in~\eqref{G-C} is in the so-called Kolmogorov
distance, which in particular implies weak convergence, that is, for every
bounded, continuous function $f$ on $\mathbb{R}$,
$\int_{\mathbb{R}} f(x) \mu _{n} (dx) \rightarrow \int_{\mathbb{R}} f(x)
\nu (dx)$.

Similar results also hold when the random variables are not
independent, but exhibit some form of weak dependence. For instance, consider
a triangular array of (dependent) random variables
$(X_{i}^{n}, i = 1,\ldots , n)_{n \in {\mathbb{N}}}$ that have a common
mean, finite variances, and exhibit \emph{asymptotic correlation decay}
in the sense that there exist positive real numbers
$\{f_{n,k}, k = 1, \ldots , n\}_{n \in {\mathbb{N}}}$ such that
$\sup_{k,n \in {\mathbb{N}}} f_{n,k} < \infty $,

\begin{equation}
\label{cordecay} \bigl\llvert \mathrm{Cov} \bigl(X^{n}_{i},X^{n}_{j}
\bigr) \bigr\rrvert \leq f_{n, \llvert i-j \rrvert } \quad \mbox{and} \quad \lim
_{n \rightarrow \infty } \frac{1}{n} \sum
_{i=1}^{n} f_{n,k} = 0,
\end{equation}
where $\mathrm{Cov}(X^{n}_{i},X^{n}_{j})$ represents the covariance of
$X^{n}_{i}$ and $X^{n}_{j}$. Then it follows from Chebyshev's inequality
\cite{Tchebichef-1867} that the normalized partial sum
$S_{n} = \frac{1}{n} \sum_{i=1}^{n} X_{i}^{n}$ satisfies
\begin{equation*}
\mathbb{P} \bigl( \bigl\llvert S_{n} - {\mathbb{E}}
\bigl[X_{1}^{1}\bigr] \bigr\rrvert > \varepsilon \bigr)
\rightarrow 0, \quad \forall \varepsilon > 0.
\end{equation*}

On the other hand, in many interesting cases one wants to analyze large
collections of strongly dependent random elements. Such an analysis is
often facilitated by  graphical model representations, which capture {\em conditional independence}
properties of the random elements via a graph.  A specific class of graphical models that will
be important for the present discussion is a Markov random field
(MRF) (a precise definition is given in Section~\ref{subs-MRF}). The theory
of MRFs and associated Gibbs measures  goes back
to the late 1960s with the pioneering works of Dobrushin
\cite{Dob68a,Dob68b} and Lanford and Ruelle \cite{LanRue69}, who were motivated
by models in statistical physics involving static  interacting random elements.
In this case a key question is efficient computation or analytical
characterization of marginal distributions of the high-dimensional ensemble of random elements.

\subsection{Questions of interest}
\label{subs-ques}

This article focuses on the dynamics of large ensembles of stochastic
processes whose interactions are governed by an underlying graph
$G = (V,E)$. Here $V$ represents a finite or countably infinite vertex
set and $E$ is a subset of unordered pairs of distinct vertices in
$V$ that represent the (undirected) edges of the graph. The graph
$G$ is always assumed to be simple (i.e., each pair in $E$ is comprised
of two distinct vertices) and locally finite, that is, for each
$v \in V$, the size of its neighborhood
$\partial _{G}(v) := \{u \in V: uv \in E\}$ is finite. The notation
$u \sim v$ will often also be used to indicate  $uv \in E$.  Given the graph $G$ and an  initial condition
$\xi = (\xi _{v})_{v \in V}$, we are interested in a collection of stochastic processes
$X^{G,\xi } = (X_{v}^{G,\xi }(t), t \geq 0)_{v \in V}$ indexed by the vertices of $G$,  that  satisfies $X_v^{G,\xi}(0) = \xi_v$ for $v \in V$, and whose   interaction structure
is governed by the graph $G$.  Specifically,   for each $v \in V$, the infinitesimal evolution of $X_v^{G,\xi}$ at any time
only depends on its own state (or history) and the states
(or histories) of neighboring particles in $G$ at that time.  Note that this includes both the  case when $X^{G,\xi}$ is Markovian, where  the infinitesimal evolution depends only on the current states of  particles, as well as non-Markovian evolutions, where the infinitesimal evolution of a particle can depend on its own history and the histories of particles in its neighbrhood.
  For conciseness, we will restrict our discussion to two types of dynamics: interacting diffusions, which are described in Section~\ref{subs-IPSdiff}, and  interacting jump processes, which are described in Section~\ref{subs-IPSjump}.
Given such (Markovian or non-Markovian) interacting processes on a large finite graph,   quantities of interest include the following:
\begin{enumerate}
\item[A.] The macroscopic behavior of the system as captured by the (global)
empirical measure process, defined by

\begin{equation}
\label{def-empmeas} \mu ^{G,\xi }(t) = \frac{1}{ \llvert V \rrvert } \sum
_{v \in V} \delta _{X^{G,\xi }_{v}(t)}, \quad t \geq 0.
\end{equation}
Note that for each $t > 0$, $\mu ^{G,\xi}(t)$ is a random probability measure
on the state space that encodes the fractions of particles taking values
in different (measurable) subsets of the state space.
\item[B.] The microscopic behavior, in particular the marginal dynamics
of a ``typical particle.'' By this we mean the dynamics of
$X_{o}^{G,\xi }$, where the vertex $o$, referred to as the root, is assumed to be chosen uniformly at
random from the finite vertex set $V$.  An important question here is to ascertain how the dynamics depends on the graph topology?
\end{enumerate}
Due to the complexity and high dimensionality of the dynamics, these quantities
are typically not amenable to exact analysis or efficient computation.
The goal instead is to identify more tractable approximations of reduced
dimension that can be rigorously justified by limit theorems, as the number of particles goes to infinity. A
desirable goal is to obtain an \emph{autonomous characterization} of the limiting
marginal dynamics of a typical particle and  evolution of the empirical
measure, which  does not refer to the full particle system dynamics.

\begin{figure}
\includegraphics{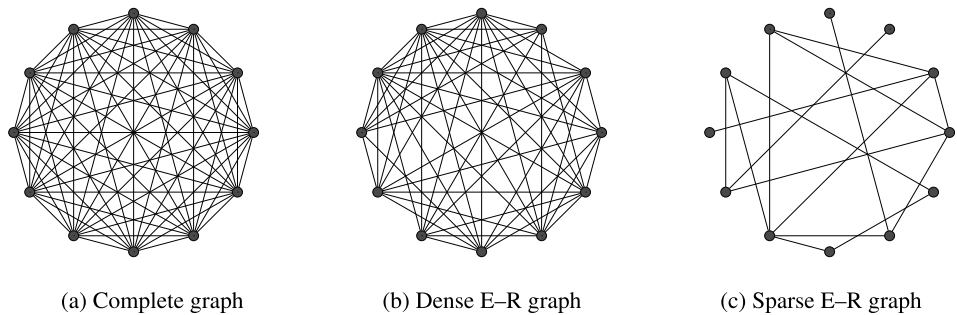}
\caption{Dense and sparse graphs} 
\label{fig-graphs}
\end{figure}

In Section~\ref{sec-stochproc} we  review the well understood case
when $G = K_{n}$, the $n$-clique or the complete graph on $n$ vertices,
in which all pairs of distinct vertices are connected by an edge; see Figure~\ref{fig-graphs}(a). For suitable initial conditions $\xi$,  convergence results for $\mu ^{K_{n},\xi }$ and
$X_{o}^{K_{n},\xi }$ in the context of interacting diffusion models go back
more than half a century to the seminal works of McKean
\cite{Mck66,Mck67}, and fall under the rubric of mean-field limits. As
briefly described in Section~\ref{sec-stochproc}, under broad conditions,
the limits, as $n \rightarrow \infty $, of both $\mu ^{K_{n},\xi }$ and
the law of $X_{o}^{K_{n},\xi }$ exist and coincide, and are described by
a certain nonlinear stochastic process. More recent work has also considered
interacting processes on certain dense random graph sequences. A generic
example of a random graph is the so-called Erd\H{o}s--R\'{e}nyi graph
$G(n,p_{n})$, which is a graph on $n$ vertices in which each pair of vertices
has an edge with probability $p_{n} \in (0,1)$ independently of all other
edges; see Figures~\ref{fig-graphs}(b) and~\ref{fig-graphs}(c) for realizations
of $G(n,p_{n})$ with $n=12$ and $p_{n}= 0.8$ and $p_{n}=0.25$, respectively.
Motivated by the study of synchronization phenomena, the work
\cite{DelGiaLuc16} considers suitably scaled pairwise interacting diffusions
on ``dense'' Erd\H{o}s--R\'{e}nyi graph $G(n,p_{n})$ sequences with divergent
average degree ($np_{n} \rightarrow \infty$), and shows that
the law of $X_{o}^{G_{n},\xi }$ converges to the same mean-field limit as
in the complete graph case. The key idea is that in this regime, particles
are only weakly interacting and become asymptotically independent, and
thus the empirical measure behaves as in the i.i.d. case~\eqref{G-C}
described in Section~\ref{subs-back}. 

The main focus of this article is on the complementary setting of interacting
stochastic processes on sequences of {\em sparse}  (possibly random) graphs, where the
(average) degrees of vertices are uniformly bounded as
$n \rightarrow \infty $. A typical example  is the Erd\H{o}s--R\'{e}nyi graph $G(n,p_{n})$ sequence when  $np_n \rightarrow c \in (0,1)$. There has been extensive
analysis of various interacting stochastic processes on deterministic sparse
graphs, originating with the work of Spitzer \cite{Spi70}, followed by
significant analysis of several Markovian models  including the      contact process, exclusion process,
and voter model. These were first studied on the $d$-dimensional lattice
(see the monographs \cite{Lig85,Lig99,KipLan99,DurLigSpiSzn12}) and then
on $d$-regular trees (e.g., \cite{Pem92,Sta01}). More recent work has also
considered processes on sparse random interaction graphs (see, e.g.,
\cite{ChaDur-contact09,DelGiaLuc16,NamNguSly19A,DurHua-contact20,GanSch20,Bhaetal21}
for an incomplete list), but none of these latter works appear to address
the main question listed above of autonomous characterization of the marginal
dynamics of a typical particle. In fact, for interacting diffusions on the sequence of sparse
Erd\H{o}s--R\'{e}nyi graphs $G_{n} = G(n,c/n)$, with
$c \in (0,\infty )$, obtaining such a characterization has remained an important open question (e.g., see \cite[p.\ 9]{DelGiaLuc16}).

The sparse regime is more challenging because particles have strong interactions,
neighboring particles remain correlated in the limit as
$n \rightarrow \infty $, and the topology of the graph has a strong influence.
Sections~\ref{sec-hydro} and~\ref{sec-marg} describe recent progress that
in particular provides a resolution of the open question in
\cite{DelGiaLuc16}. The article concludes in Section~\ref{sec-open} with
generalizations and open questions. The work on both mean-field models
and various other aspects of interacting particle systems is so extensive
that it will be impossible to be representative in this short
article. Instead, I hope to just provide enough pointers for the reader
to get a flavor of the classical results and set the context for  more
recent results.  Monographs covering
various aspects of interacting particles systems include
\cite{Lig85,Lig99,Szn91,BiaDur93,KipLan99,Kol10,Geo11,dembo-montanari}.

\section{Classical mean-field results for interacting stochastic processes}
\label{sec-stochproc}

Given a (simple, locally finite, undirected) graph $G = (V,E)$ and
$v \in V$, we use
$\mathrm{cl}_{G}(v) := \{v\} \cup \partial _{G}(v)$ to denote the closure
of $v$ in $G$. Note that
$|\mathrm{cl}_{G}(v)|$ is always finite, where $|A|$ denotes the cardinality
of a set $A$. We now describe the dynamics of locally interacting diffusions
and interacting jump processes. Rather than provide the most general setting, we make  simplifying assumptions whenever
convenient to illustrate the key issues.

\subsection{Interacting diffusions}
\label{subs-IPSdiff}

Given an initial condition
$\xi = (\xi _{v})_{v \in V} \in \mathbb{R}^{V}$ with
${\mathbb{E}}[\xi _{v}^{2}]< \infty $ for every $v \in V$, consider the
collection $X^{G,\xi } = \{X_{v}^{G,\xi }\}_{v \in V}$ of diffusive particles,
indexed by the nodes of the graph~$G$,  that evolve according to the following
coupled system of stochastic differential equations (SDEs):
\begin{equation}
\label{SDE-diff2} dX_{v}^{G,\xi } (t) = b
\bigl(t,X_{v}^{G,\xi }(t), \mu _{v}^{G,\xi }(t)
\bigr)\, dt + dW_{v}(t), \quad X_{v}^{G,\xi } (0)
= \xi _{v}, \quad t > 0, v \in V,\quad
\end{equation}
where $(W_{v})_{v \in V}$ are i.i.d. standard Brownian motions independent
of $(\xi _{v})_{v\in V}$, and for any vertex that is not isolated, $\mu _{v}^{G,\xi }(t)$ represents the \emph{local
empirical measure} of a neighborhood of $v$ at time $t \geq 0$,
\begin{equation*}
\mu _{v}^{G,\xi }(t) = \frac{1}{ \llvert \partial _{G}(v) \rrvert } \sum
_{u \in
\partial _{v}} \delta _{X_{u}^{G,\xi }(t)},
\end{equation*}
and $b$ is a drift coefficient that is sufficiently regular to ensure  that the SDE~\eqref{SDE-diff2} has a unique weak solution. (When $v$ is  isolated,   the precise definition of $\mu_v^{G,\xi}$ is not so important; it can be set  equal to an arbitrary quantity.)

A special case of interest is when $b$ has  linear dependence on the measure
term, say, of the form

\begin{equation}
\label{pairwise} b(t,x,\nu ) = \int_{\mathbb{R}} \beta (t,x,y) \nu
(dy),
\end{equation}
for some interaction potential
$\beta : \mathbb{R}_{+} \times \mathbb{R}^{2} \to \mathbb{R}$ that is symmetric
in the last two variables. In this case, system~\eqref{SDE-diff2} reduces
to the following system of pairwise interacting diffusions:
\begin{equation}
\label{SDE-diff} dX_{v}^{G,\xi } (t) = \frac{1}{ \llvert \partial _{G}(v) \rrvert }
\sum_{u \sim v} \beta \bigl(t,X_{v}^{G,\xi }(t),
X_{u}^{G,\xi }(t)\bigr) dt + dW_{v}(t), \quad v
\in V, t > 0,
\end{equation}
which models phenomena in different fields, including statistical physics
and neuroscience \cite{Der03,RedRoeRus10,LucSta14}. The trajectories of
each particle lie in the space $\mathcal{C}$ of continuous real-valued
functions on $[0,\infty )$, which we endow with the topology of uniform
convergence on compact sets.

\subsection{Mean-field limits and nonlinear diffusion processes} 
\label{subs-MFdiff}

Now consider the SDE~\eqref{SDE-diff} with $G = K_{n}$, the complete graph,
and assume without loss of generality that $G$ has vertex set
$\{1, \ldots , n\}$. We present a sufficient condition on the drift
under which one can establish a standard mean-field result. Given any
$p \geq 1$ and Polish space $E$, the Wasserstein-$p$ metric on $E$ is defined
as follows:
\begin{align}
\label{Wassp} \mathcal{W}_{E,p}(\nu ,\tilde{\nu }) := \inf
_{\pi } \biggl(\int_{E
\times E}
d^{p}(x,y)\pi (dx,dy) \biggr)^{1/p},
\end{align}
where the infimum is over all couplings $\pi $ of $\nu $ and
$\tilde{\nu }$, namely probability measures $\pi $ on $E^{2}$ with first
and second marginals $\nu $ and $\tilde{\nu }$, respectively. Let
$\mathcal{P}^{p} (E)$ be the space of probability measures on $E$ equipped
with the Wasserstein-$p$ metric $\mathcal{W}_{E,p}$.

\begin{assumption}
\label{ass-MFdiff}
Suppose that $b$ is bounded and for every $t > 0$, the map
$\mathbb{R}\times \mathcal{P}^{2}(\mathbb{R}) \ni (x,\nu ) \mapsto b(t,x,
\nu ) \in \mathbb{R}$ is Lipschitz continuous,  uniformly with respect to $t$ in compact
subsets of  $\mathbb{R}_{+}$.
\end{assumption}

Note that Assumption~\ref{ass-MFdiff} is satisfied
when the drift $b$ is of the form~\eqref{pairwise}, where the interaction potential  $\beta$  is such that
$\mathbb{R}^2 \ni (x,y) \to \beta(t,x,y)$  is  Lipschitz continuous and bounded,   uniformly with respect to $t$ in compact subsets
of $\mathbb{R}_{+}$.

\begin{theorem}
\label{th-MFdiff}
Suppose Assumption~\ref{ass-MFdiff} holds, and there exists
$\mu _{o} \in \mathcal{P}^{2}(\mathbb{R})$ such that the initial conditions
$(\xi ^{n}_{i})_{i=1, \ldots , n}$, $n \in {\mathbb{N}}$, satisfy
\begin{equation}
\label{ic-poc} {\mathbb{E}} \Biggl[\mathcal{W}_{\mathbb{R},2} \Biggl(
\frac{1}{n} \sum_{i =1}^{n}
\delta _{\xi _{i}^{n}} , \mu _{o} \Biggr) \Biggr] \rightarrow 0
\quad \mbox{as } n \rightarrow \infty .
\end{equation}
Then there is a unique strong solution to the SDE
\begin{equation}
\label{diff-nonlinear} d X_{o} (t) = b\bigl(t,B(t), \mu (t)\bigr)\, dt +
dB(t), \quad \mu (t) = \mathscr{L}\bigl(X_{o}(t)\bigr), \quad t > 0,
\end{equation}
with $\mathscr{L}(X_{o}(0)) = \mu _{0}$. Moreover, if for each
$n \in {\mathbb{N}}$, $X^{n} := X^{K_{n},\xi ^{n}}$ is the unique solution
to the SDE~\eqref{SDE-diff2}, then the global empirical measure
$\mu ^{n} := \mu ^{K^{n},\xi ^{n}}$ defined in~\eqref{def-empmeas} satisfies
\begin{equation}
\label{timet-conv} \lim_{n \rightarrow \infty } {\mathbb{E}} \Bigl[ \sup
_{s \in [0,t]} \mathcal{W}_{\mathbb{R},2} \bigl(\mu
^{n}(s), \mu (s)\bigr) \Bigr] = 0, \quad \forall t > 0.
\end{equation}
Furthermore, for any $k \in \mathbb{N}$ and $t > 0$, the law of
$(X_{1}^{n}(t), \ldots, X_{k}^{n}(t))$ converges weakly to the product
$(\mu (t))^{\otimes k}$, that is, for all bounded continuous functions
$f_{i}: \mathbb{R}\to \mathbb{R}$, $i = 1, \ldots , k$,
\begin{equation}
\label{timet-chaos} \lim_{n \rightarrow \infty } {\mathbb{E}}\bigl[
f_{1}\bigl(X_{1}^{n}(t)\bigr) \ldots
f_{k}\bigl(X_{k}^{n}(t)\bigr) \bigr] =
\prod_{i=1}^{k} \int
_{\mathbb{R}} f_{i} (x) \mu (t) (dx).
\end{equation}
\end{theorem}

If there were no interaction, $b \equiv 0$, then the theorem would simply
be a (functional) strong law of large numbers result. However, even when
$b \not\equiv 0$, the particles are only weakly interacting because the symmetry
of the interaction ensures that the influence of any particle on the drift
of another particle is $O(1/n)$, which vanishes in the limit. The property~\eqref{timet-chaos} that any finite subset of random variables from
$\{X^{n}_{i}(t), i = 1, \ldots , n\}_{n \in \mathbb{N}}$ are asymptotically
independent is referred to as \emph{chaoticity}, and is well known to be
equivalent to the convergence of $\mu ^{n}(t)$ to a deterministic law
\cite[Proposition 2.2]{Szn91}. Now,~\eqref{ic-poc} implies that the initial
conditions are chaotic.  Thus Theorem~\ref{th-MFdiff} asserts that the dynamics are such that this chaoticity also holds for positive times $t > 0$, a phenomenon referred to as \emph{propagation of chaos}.  In turn,  this   leads to an \emph{autonomous} description
of the limiting marginal process $X_{o}$, which is   a Markov process whose infinitesimal evolution at any  time $t$ also depends on its own law $\mu (t)$ at that time.  As a result, the
forward Kolmogorov equation  (or master equation), which is the
partial differential equation (PDE) describing the evolution
of the marginal law $\mu$,  is nonlinear.  Consequently, such a process
is referred to as a \emph{nonlinear} Markov process. When the drift has
the form~\eqref{pairwise}, under suitable conditions it can be shown that the law $\mu (t)$ is absolutely continuous with respect to Lebesgue measure and that its density satisfies the granular media equation
\cite{Mck66}. Thus, PDE techniques can be useful for studying nonlinear Markov processes (see, e.g., \cite{BarRoc20}).

There are many different approaches to establishing mean-field limits,
including PDE analysis, fixed point arguments, martingale techniques and stochastic coupling constructions.  First, one needs to  establishing well-posedness of the
nonlinear SDE~\eqref{diff-nonlinear}.  An analytical approach to this problem  entails  proving  uniqueness of the nonlinear PDE describing the evolution of the marginal law.
Another, more probabilistic,  approach is to first consider the mapping that
takes any continuous measure flow
$t \mapsto \nu (t) \in \mathcal{P}^{2}(\mathbb{R})$ to the measure flow
$t \mapsto \mathscr{L}(X^{\nu }(t))$, where $X^{\nu }$ is the unique solution to
the SDE in~\eqref{diff-nonlinear} when $\mu $ is replaced with
$\nu $.   Observing that the flow $t \mapsto \mathscr{L}(X_{o}(t))$ must be
a fixed point for this mapping,  well-posedness is equivalent to uniqueness of the fixed point of this mapping.   The latter can be established by showing the mapping is a contraction by exploiting the Lipschitz continuity of the drift.   Given well-posedness, the coupling approach to proving convergence proceeds by first defining $\bar{X}^{n}$ to be the $n$-dimensional process whose every
coordinate is an independent copy of the nonlinear
process {$X_{o}$.}  Then one couples this process with the original process $X^{n}$ so that they are both driven by the same Brownian motions.
Using
It\^{o}'s formula,  the Lipschitz condition on the drift and standard estimates, one can then  show that the $\mathcal{W}_{\mathbb{R},2}$ distance between the empirical measures
of $X^n$ and $\bar{X}^n$ vanishes as $n \rightarrow \infty.$  Since the strong law of large numbers ensures that the empirical measure of the latter converges to the law of $X_{o}$, which is equal to $\mu$, this concludes the proof.  An alternative approach to proving convergence is to first  use the generator of the Markov process $X^{n}$ to identify martingales involving
the empirical measure process $\mu ^{n}$, next show that the sequence $\{\mu ^{n}\}$ is relatively compact (or tight), then characterize any subsequential limit satisfies
what is known as a nonlinear martingale problem, and finally establish
well-posedness of the latter \cite{oelschlager1984martingale,Gar88}.

\begin{remark}
One can consider more general dynamics where both the drift and diffusion
coefficients are functions of the current state and the empirical measure
process, as well as non-Markovian versions that depend on the history of
the process.
\end{remark}

\subsection{Interacting jump processes and their mean-field limits} 
\label{subs-IPSjump}

\subsubsection{Description of dynamics}
\label{188.sec2.3.1}

We will also be interested in interacting pure jump processes, which describe
models in statistical physics, engineering, epidemiology and
the dynamics of opinion formation \cite{Lig99,BiaDur93}. For
concreteness, consider the voter model \cite{Lig99} that aims to capture opinion
dynamics, in which each particle takes values in the state space
$\mathcal{X} = \{0,1\}$ that represents two possible opinions. The allowed
transitions or jump directions of a particle lie in the set
$\mathcal{J} = \{1,-1\}$.  The rate at which any  particle   changes its opinion is equal to the fraction of its neighbors with the opposite opinion.  Note that the dependence of the rate on the neighboring states is symmetric.  More generally, when the state of the system is $(x_v)_{v\in V}$, the jump rate  of a particle at $v$  could be a more complicated symmetric functional of the neighboring states
$(x_u)_{u \sim v}$ and also depend on time $t$, in addition to its own state $x_v$.
This symmetric dependence on neighboring states is most  succinctly  captured by saying the rate is a functional of the {\em unnormalized} empirical measure $\theta_v = \sum_{u \sim v} \delta_{x_u}$ of the neighboring states.  Note that $\theta_v$ lies in the space
$\mathcal{M}(\mathcal{X})$
of locally finite nonnegative integer-valued measures on $\mathcal{X}$.

In the general setup, we consider a finite state space
$\mathcal{X}$, a  subset
$\mathcal{J} \subset \{i-j: i, j \in \mathcal{X}\}$ of possible jump directions, and  a collection of jump rate functions
$\bar{r}_{j}: \mathbb{R}_+ \times \mathcal{X}\times \mathcal{M}(\mathcal{X}) \to
\mathbb{R}_{+}$, $j \in \mathcal{J}$.
Given a (simple) finite graph
$G = (V,E)$ and  initial condition $\xi = (\xi_v)_{v \in V} \in {\mathcal X}^V$, the $\mathcal{X}^{V}$-valued process representing the configuration
of the associated IPS evolves according to the following system of (jump) SDEs:

\begin{equation}
\label{SDE-jump} X^{G,\xi }_{v}(t) = \xi _{v} +
\sum_{j \in \mathcal{J}} j \int_{(0,t]
\times \mathbb{R}_{+}}
\mathbb{I}_{\{r\leq \bar{r}_{j}(s,X^{G,\xi }_{v}(s-),
\theta _{v}^{G,\xi }(s-))\}} N_{v}(ds,dr),
\quad t \geq 0, v \in V,
\end{equation}
where $(N_{v})_{v \in V}$, are i.i.d. Poisson random measures on
$\mathbb{R}_{+}^{2}$ with intensity measure $\Leb ^{2}$, where
$\Leb$ represents Lebesgue measure on $\mathbb{R}$, and for each
$s \geq 0$,
$\theta^{G,\xi}_{v}(s)$ is the random (unnormalized) empirical measure corresponding to the states of the particles in the neighborhood
of $v$ at time $s$:
\begin{equation}
\label{theta_v} \theta _{v}^{G,\xi }(s) := \sum
_{u \sim v} \delta _{X^{G,\xi }_{u}(s)}, \quad v \in V, s \geq 0.
\end{equation}
The SDE~\eqref{SDE-jump} captures a simple evolution.  For any $j \in {\mathcal J}$ and time
$t$, the particle at a node $v$  makes a transition from its state $X^{G,\xi }_{v}(t-)$ to $X^{G,\xi }_{v}(t-) + j$  at
a rate $\bar{r}_{j}(t,x, \theta _{v}^{G,\xi }(t-))$ that
depends on the current time, the state of the node just prior to the current time, and symmetrically
on the states of neighboring nodes just prior to the current time, as encoded by
$\theta _{v}^{G,\xi }(t-)$. 
Use of the unnormalized measure
$\theta _{v}^{G,\xi }(t)$ instead of the empirical measure allows one to capture
a broader class of models in which jump rates depend on the number of neighboring
nodes in particular states (and not just their fractions), as is the case for models like the contact process \cite{Lig99}.
Note that the
trajectory of each particle lies in the c\`{a}dl\`{a}g space
$\mathcal{D}$ of right continuous ${\mathcal X}$-valued functions on
$[0,\infty )$ that have finite left limits on $(0,\infty )$.

The solution $X^{G,\xi }$ to the jump SDE~\eqref{SDE-jump} is a Markov jump
process and so its law can also be characterized via the associated
{\em infinitesimal generator} \cite{Lig85}: for functions
$f: \mathcal{X}^{V} \mapsto \mathbb{R}$,
\begin{align*}
\mathcal{A}_{t} f(x) &= \lim_{h \downarrow 0}
\frac{\mathbb{E}[ f(X_{t+h}^{G,\xi }) - f(X_{t}^{G,\xi })|X_{t}^{G,\xi } = x]}{h}
\\
&= \sum_{j \in \mathcal{J}, v \in V} \bar{r}_{j}
\biggl(t,x_{v},\sum_{u
\sim v} \delta
_{x_{u}} \biggr) \bigl[f (x + j e_{v} ) - f(x) \bigr],
\quad t > 0, x \in \mathcal{X}^{V},
\end{align*}
where $e_{v} \in \{0,1\}^{V}$ is the vector with $1$ in the $v$th coordinate
and $0$ elsewhere. However, the jump SDE representation
in~\eqref{SDE-jump} is more convenient for  generalizations to non-Markovian processes (see \cite{GanRam-Hydro22}).
 Furthermore, the jump SDE formulation is also better suited to describing the form of limiting marginal dynamics on sparse graphs, as
described in Section~\ref{subs-margjump}.

\subsubsection{Mean-field limits and nonlinear jump processes} 
\label{subs-hydrojump}

Mean-field results analogous to Theorem~\ref{th-MFdiff} also hold in the
jump setting
under the following regularity   assumption on the jump rate functions:

\begin{assumption}
\label{ass-MFjump2}
For each $j \in \mathcal{J}$, the jump rate function takes the form
$\tilde{r}_{j}(t,x,\theta ) = \hat{r}_j(t,x, \frac{1}{\theta (\mathcal{X})}
\theta )$ when $\theta({\mathcal X}) \neq 0$, and $\tilde{r}_{j}(t,x,\theta ) = 0$ otherwise,  where the function
$\hat{r}_{j}: \mathbb{R} \times \mathcal{X}\times \mathcal{P}^{1}(
\mathcal{X}) \mapsto \mathbb{R}_{+}$ is  such that ${\mathcal P}^1({\mathcal X}) \ni \nu \mapsto \hat{r}(t,x,\nu)$ is Lipschitz continuous, uniformly for $x \in {\mathcal X}$ and $t$ in compact subsets
of $\mathbb{R}_{+}$.
\end{assumption}

Assumption \ref{ass-MFjump2} reflects the fact that in the mean-field setting, the dependence of the jump rates on the neighboring particles must be a sufficiently regular function of the usual (normalized) empirical measure.
The following result is established in
\cite[Theorem 2]{oelschlager1984martingale}; see also \cite{Kol10}.

\begin{theorem}
\label{th-MFjump}
Suppose Assumption~\ref{ass-MFjump2} holds and  the initial conditions
are chaotic, that is,
$\frac{1}{n} \sum_{i=1}^{n} \delta _{\xi _{i}^{n}}$ converges in the total
variation metric to a deterministic limit $\mu _{0}$, then
$\mu ^{K_{n},\xi ^{n}}(t)$ converges weakly to
$\mathscr{L}(X_{o}(t))$ where
$X_{o}(t) = X^{\mu _{0}}(t)$, $t \geq 0$, is the unique solution to the following
nonlinear jump SDE: $\mathscr{L}(\mathcal{X}_{o}(0)) = \mu _{0}$, and for
$t \geq 0$,
\begin{align}
\label{jump-nonlinear} X_{o}(t) & = X_{o}(0) + \sum
_{j \in \mathcal{J}} j \int_{(0,t]
\times \mathbb{R}_{+}}\mathbb{I}_{\{r\leq \bar{r}_{j}(s,X_{o}(s-),
\mu (s-))\}}
N(ds,dr),
\\
\nonumber
\mu (t) & = \mathscr{L}\bigl(X_{o}(t)\bigr),
\end{align}
where $N$ is a Poisson process on $\mathbb{R}^2_+$ with intensity
$\Leb ^{2}$, independent of $X_{o}(0)$. Furthermore, for any
$k \in \mathbb{N}$, the law of
$(X_{1}^{K_{n},\xi }, \ldots , X_{k}^{K_{n},\xi })$ on
$\mathcal{D}^{k}$ converges weakly to
$(\mathcal{L}(X_{o}))^{\otimes _{k}}$.
\end{theorem}

Just as the evolution of the  law of the mean-field diffusion limit in~\eqref{diff-nonlinear} can be characterized by a nonlinear PDE, the evolution
of the law of the nonlinear jump process $X_{o}$ in~\eqref{jump-nonlinear} can be characterized as the unique solution to its forward equation, which is now a
nonlinear integrodifferential equation.

\begin{remark}
Theorems~\ref{th-MFdiff} and~\ref{th-MFjump} are meant to only provide
a flavor of mean-field results. While a survey of mean-field limits is
not the current focus, it is worth mentioning that in both
the diffusive and jump process settings, one can obtain mean-field limits
under weaker assumptions and for much more general dynamics where the diffusion
coefficient is also a function of the current state and empirical
measure process, as well as non-Markovian versions where the drift coefficient
or jump rates depend on the history of the process (see, e.g.,
\cite{MehSchStaZag20} for propagation of chaos results on interacting non-Markovian
jump diffusions and \cite{BalPerRei22} for a large deviations analysis
of non-Markovian weakly interacting diffusions).
\end{remark}

\begin{figure}[b!]
\vspace*{-12pt}
\includegraphics{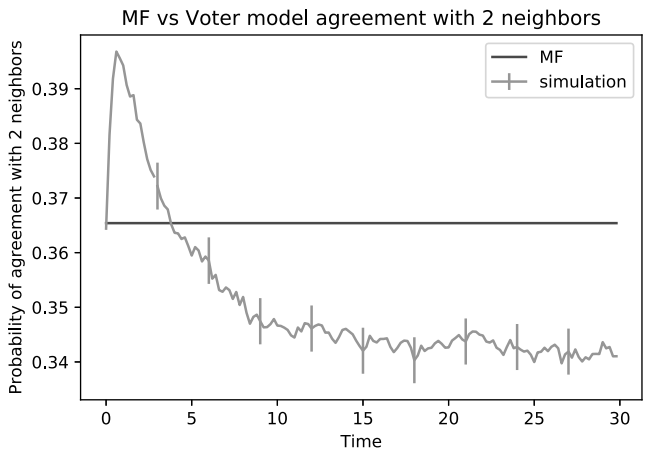}
\caption{Simulations and mean-field (MF) approximations for the voter model
on a $3$-tree.}
\label{fig-votermodel}
\vspace*{-6pt}
\end{figure}

\vspace*{-3pt}

\subsection{Limitations of mean-field approximations}
\label{188.sec2.4}

The mean-field limit theorems established in Theorems~\ref{th-MFdiff} and~\ref{th-MFjump} indicate that the law of the nonlinear Markov processes
$X_{o}$ in~\eqref{diff-nonlinear} and~\eqref{jump-nonlinear}, respectively, can be
used to approximate quantities of interest for interacting diffusions or
jump processes on finite graphs. In particular, consider the voter model described in
Section \ref{188.sec2.3.1}. Its jump rates take the explicit form
  \begin{align*}
 \bar{r}_{1}(t,x, \theta ) & = \frac{\mathbb{I}_{\{x = 0\}}}{ \llvert \theta (\mathcal{X}) \rrvert } \int
 _{
 \mathcal{X}} y \theta (dy), \quad \bar{r}_{-1}(t,x,
 \theta ) = \frac{\mathbb{I}_{\{x = 1\}}}{ \llvert \theta (\mathcal{X}) \rrvert }{\int_{
 \mathcal{X}} (1-y) \theta
 (dy)}.
 \end{align*}
In this case,
one would expect that the dynamics of
$X_{o}$ in~\eqref{jump-nonlinear} with these reates could provide an   approximation for the probability
of agreement of any two neighboring particles in the voter model on a sufficiently
large complete graph.  However, for lack
of a better alternative,  mean-field approximations are used even for dynamics on other graphs.  While these approximations may do reasonably well on
dense graphs (where vertices have high degrees) \cite{DelGiaLuc16}, they
can be very inaccurate on sparse graphs. Figure~\ref{fig-votermodel} plots
the evolution of the probability that the state of the root agrees with  precisely two of its neighbors for the voter model on a
rooted $3$-regular tree with $9$ generations, given  at time zero, each particle
independently has an opinion $1$ with probability $0.3$. The vertical bars in Figure~\ref{fig-votermodel} provide confidence intervals for the simulation. The mean-field approximation
assumes neighboring vertices are independent, and thus performs poorly.
Ad hoc refinements of the mean-field approximation that take into account
correlations also remain inaccurate in this setting. This strongly motivates
the development of a convergence theory for the empirical distribution
and marginal dynamics on sparse graph sequences that could lead to more
principled approximations.

\section{Interacting processes on sparse graphs: hydrodynamic limits} 
\label{sec-hydro}

We now turn to interacting processes
$(X^{G_{n},\xi ^{n}})_{n \in \mathbb{N}}$ on sparse graph sequences
$(G_{n})_{n \in \mathbb{N}}$ with initial conditions
$(\xi ^{n})_{n \in \mathbb{N}}$. Assume each $G_{n}$ is finite and
$o_{n}$ is a vertex chosen uniformly at random from the vertices of
$G_{n}$. Unlike in the case of the complete graph (or even dense graph
sequences), the degree of a vertex remains bounded and so neighboring vertices do not become asymptotically independent.  Thus, the number of neighbors becomes important and so it is clear that one cannot expect $(X^{G_{n},\xi ^{n}}_{o_{n}})_{n \in \mathbb{N}}$ to have a limit just
by sending the number of vertices $n$ to infinity, without imposing any
additional consistency requirements on the graphs in the sequence. This
leads to the following questions:
\begin{enumerate}
\item[Q1.] For what graph sequences $(G_{n})_{n \in {\mathbb{N}}}$ would
one expect $(X^{G_{n},\xi ^{n}}_{o_{n}})_{n \in {\mathbb{N}}}$ to have a
limit?
\item[Q2.] For such sequences, will
$(\mu ^{G_{n},\xi ^{n}})_{n \in {\mathbb{N}}}$ converge to a deterministic
limit?
\item[Q3.] When $(\mu ^{G_{n},\xi ^{n}})_{n \in {\mathbb{N}}}$ converges
to a deterministic limit, will this limit always coincide with the limit
law of $X^{G_{n},\xi ^{n}}_{o_{n}}$?
\item[Q4.] Is there an autonomous reduced-dimension description of the
limit of the marginal $X^{G_{n},\xi ^{n}}_{o_{n}}$ whenever this limit
exists?
\end{enumerate}

In light of the first question above, we review a natural notion of convergence of sparse graphs
called local convergence in Section~\ref{subs-locconv}. This notion was
used to study asymptotic properties of static models (Gibbs measures) of
discrete-valued marked random graphs in \cite{dembo-montanari}.

\subsection{Local convergence of sparse graph sequences}
\label{subs-locconv}

Given a graph $G=(V,E)$ and two vertices $u, v \in V$, a path of length
$n$ between $u$ and $v$ is a sequence
$u=u_{0},u_{1}, \ldots , u_{n}=v$ such that $u_{i-1} \sim u_{i}$ for every
$i = 1, \ldots , n$. A graph is said to be connected if there exists a
finite path between any two vertices and the graph distance between two
vertices is the minimum length of a path between them. A~rooted graph $(G,o)$ is a graph $G = (V,E)$ with a special vertex $o \in V$, referred to as the root.  A useful notion
of convergent sequences of (connected) rooted sparse graphs is that of \emph{local convergence}, which
 was introduced by Benjamini and Schramm
\cite{benjamini2001recurrence}. Other references on local convergence include
\cite{Bordenave2016,aldous-steele}.  We first introduce some terminology that is required to define local convergence. An \emph{isomorphism} from one rooted
graph $(G_{1},o_{1})$ to another $(G_{2},o_{2})$ is a bijection
$\varphi $ from the vertex set of $G_{1}$ to that of $G_{2}$ such that
$\varphi (o_{1})=o_{2}$ and such that $(u,v)$ is an edge in $G_{1}$ if
and only if $(\varphi (u),\varphi (v))$ is an edge in $G_{2}$. Two rooted
graphs are said to be \emph{isomorphic} if there exists an isomorphism between them.
Let $\mathcal{G}_{*}$ denote the set of isomorphism classes of connected
rooted graphs.  We will also need to consider convergence of graphs that  carry "marks" representing the initial condition or trajectory of the state dynamics at that vertex.
With that in mind, given a Polish space $\mathcal{S}$, we define a
\emph{$\mathcal{S}$-marked rooted graph} to  be a tuple $(G,x,o)$, where
$(G,o)$ is a rooted graph and
$x=(x_{v})_{v \in G} \in \mathcal{S}^{G}$ is a vector of marks, indexed
by the vertices of $G$. We say that two marked rooted graphs
$(G_{1},x^{1},o_{1})$ and $(G_{2},x^{2},o_{2})$ are \emph{isomorphic} if
there exists an isomorphism $\varphi $ from the rooted graph
$(G_{1},o_{1})$ to the rooted graph $(G_{2},o_{2})$ that maps the marks of $(G_{1}, o_{1})$ to the marks
of $(G_{2},o_{2})$ (i.e., for which $x^{2}_{\varphi (v)} = x^{1}_{v}$ for all
$v \in G$). Let $\mathcal{G}_{*}[\mathcal{S}]$ denote the set of isomorphism
classes of $\mathcal{S}$-marked rooted graphs.

We  now define the  topologies of local convergence on the spaces $\mathcal{G}_{*}$ and $\mathcal{G}_{*}[\mathcal{S}]$. For $r \in {\mathbb{N}}$ and $(G,o)\in \mathcal{G}_{*}$, let $B_{r}(G,o)$ denote the induced subgraph
of $G$ (rooted at $o$) containing  those vertices with (graph) distance
at most $r$ from the root $o$. The distance between $(G_{1},o_{1})$ and
$(G_{2},o_{2})$ in $\mathcal{G}_*$ is defined to be $1/(1+\bar{r})$, where
$\bar{r}$ is the supremum over $r \in {\mathbb{N}}_{0}$ such that
$B_{r}(G_{1},o_{1})$ and $B_{r}(G_{2},o_{2})$ are isomorphic, where we
interpret $B_{0}(G_{i},o_{i}) = \{o_{i}\}$.   Now, let $d$ denote a metric  that induces the Polish topology on $\mathcal{S}$. We then metrize $\mathcal{G}_{*}[\mathcal{S}]$  by similarly defining the  distance between two $\mathcal{S}$-marked
graphs $(G_{i},x^{i},o_{i})$, $i = 1, 2$, to be
$1/(1+\bar{r})$, where now $\bar{r}$ is the supremum over
$r \in {\mathbb{N}}_{0}$ such that there exists an isomorphism
$\varphi $ from $B_{r}(G_{1},o_{1})$ to $B_{r}(G_{2},o_{2})$ for which
$d(x^{1}_{v},x^{2}_{\varphi (v)}) \le 1/r$ for all
$v \in B_{r}(G_{1},o_{1})$. Under the respective topologies,
$\mathcal{G}_{*}$ and $\mathcal{G}_{*}[\mathcal{S}]$ are Polish spaces
(see \cite[Lemma 3.4]{Bordenave2016} or
\cite[Appendix~A]{LacRamWuLWC21}). For any Polish space
$\mathcal{S}$, let $C_{b}(\mathcal{S})$ denote the space of bounded
continuous functions on~$S$.\looseness=-1

We will always assume the spaces $\mathcal{G}_{*}$ and
$\mathcal{G}_{*}[\mathcal{S}]$ are equipped with their Borel $\sigma $-algebras. One can then talk about weak convergence or convergence in distribution of random graphs and random marked graphs as random elements in $\mathcal{G}_*$ or $\mathcal{G}_*[\mathcal{S}].$   Specifically, a sequence of random $\mathcal{G}_*$-valued random elements $\{(G_n,o_n)\}$ is said to converge in distribution in the local weak sense to a $\mathcal{G}_*$-valued limit $(G,o)$ if for every bounded continuous function $f:\mathcal{G}_* \rightarrow \mathbb{R}$, $\mathbb{E}[f(G_n,o_n)] \rightarrow \mathbb{E}[f(G,o)].$  Likewise,  convergence in distribution in the local weak  sense of (isomorphism classes of) random $\mathcal{S}$-marked graphs is equivalent to weak convergence on the space $\mathcal{G}_*[\mathcal{S}]$.

 \begin{figure}
   \hspace{-0.6in}
   (a) \hspace{-0.25in}
 \begin{minipage}{0.5\linewidth}
 \tikzstyle{every node}=[fill,scale=0.5]
\begin{center}
{\ } \\ \vspace{0.7em}
\begin{tikzpicture}
\def \n {12}
\def \radius {1.7cm}
\def \margin {1} % margin in angles, depends on the radius
  \node[label=right:{$o$},draw, circle] at ({0}:\radius) {};
\foreach \s in {1,...,\n}
{
  \node[draw, circle] at ({360/\n * (\s - 1)}:\radius) { };
    \draw ({360/\n * (\s - 1)+\margin}:\radius)
    edge ({360/\n * (\s)-\margin}:\radius);
}
\end{tikzpicture}
\end{center}
   \end{minipage}%
   \begin{minipage}{0.1\linewidth}
  {\ } \\  {\ } \\  \vspace{1.3em}
\hspace{-0.25in} $\longrightarrow$ \hspace{0.25in}
   \end{minipage}
   \hspace{0.3in}
      \begin{minipage}{0.2\linewidth}
\begin{center}
$\vdots \ \ \ $ \\ \vspace{0.7em}
\tikzstyle{every node}=[fill,scale=0.5]
  \begin{tikzpicture}[scale=0.5]
    \node[label=right:{$o$},draw, circle] at (0,4) {};
    \foreach \i in {1,...,4}
      \foreach \j in {1,...,7}{
        \node[draw, circle] at (0,\j) {};
        }
    \foreach \i in {1,...,1}
      \foreach \j in {1,...,6}{
        \draw (0,\j) edge (0, \j+1);
      }
  \end{tikzpicture}
{\ } \\  $\vdots \ \ \ $
\end{center}
   \end{minipage}\\
   \hspace{-0.15in}
(b)
      \begin{minipage}{0.4\linewidth}
\begin{center}
\tikzstyle{every node}=[fill,scale=0.5]
\begin{tikzpicture}
\def \n {16}
\def \p {41}
\def \radius {1.7cm}
\def \margin {1} % margin in angles, depends on the radius
  \node[label=right:{$o$},draw, circle] at ({0}:\radius) {};
\foreach \s in {1,...,\n}
{
  \node[draw, circle] at ({360/\n * (\s - 1)}:\radius) { };
}
\foreach \s in {1,...,\n}
{
  \draw ({360/\n * (mod(206*\s,\p) - 1)+\margin}:\radius)
    edge ({360/\n * (mod(207*\s,\p)-\margin}:\radius);
}
\foreach \s in {1,...,\n}
{
  \draw ({360/\n * (mod(101*\s,\p) - 1)+\margin}:\radius)
    edge ({360/\n * (mod(102*\s,\p)-\margin}:\radius);
}
\end{tikzpicture}
\end{center}
\end{minipage}
   \begin{minipage}{0.1\linewidth}
   {\ } \\  {\ } \\  \vspace{1.3em}
 $\longrightarrow$
   \end{minipage}
\begin{minipage}{0.4\linewidth}
\begin{center}
\tikzstyle{every node}=[fill,scale=0.5]
\begin{tikzpicture}[circle,level distance=0.8cm,
  level 1/.style={sibling distance=1.5cm},
  level 2/.style={sibling distance=1cm},
  level 3/.style={sibling distance=0.5cm}]
  \node[scale=0.3,label=above:{$o$}] {root}
        child {node {}
      child {node {}}
      child {node {}
        child {node {}}
      }
      child {node {}
        child {node {}}
        child {node {}}
        child {node {}}
      }
    }
    child {node {} }
    child {node {}
      child {node {}}
      child {node {}
                child {node {}}
                child {node {}}
                child {node {}}
                child {node {}}
      }
    };
\end{tikzpicture}
$\hphantom{x} \quad\quad\quad\quad \ \ \ \vdots \quad\quad\quad\quad\quad\quad \ \ \vdots$
\end{center}
\end{minipage}
\caption{Local convergence: (a) cycle  to  infinite line; (b) Erd\H{o}s--R\'enyi  graph to  a UGW tree.}
\label{fig-locconv}
\end{figure}
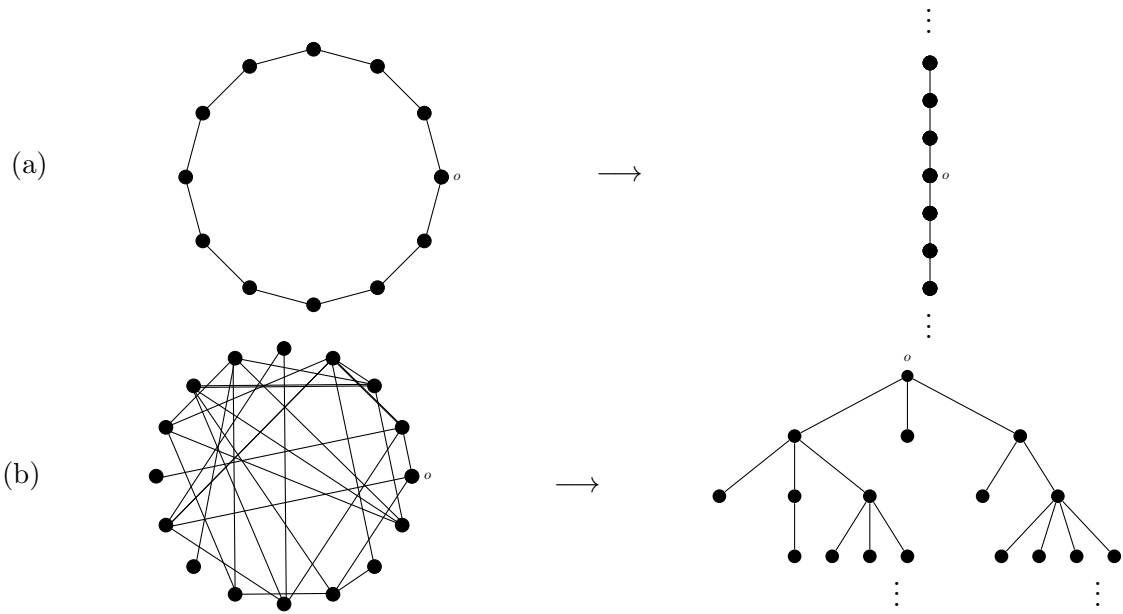

\begin{remark}
\label{rem-locconv}
Figure~\ref{fig-locconv} illustrates two generic examples of locally convergent
graph se\-quences. Let $G_{n}$ be the $n$-cycle, which is the connected graph
on $n$ vertices where every vertex has degree $2$, along with the root
$o_{n}$ chosen uniformly at random from the $n$ vertices. Then
$(G_{n},o_{n})$ converges weakly in $\mathcal{G}_{*}$ to a  infinite line graph rooted at some fixed vertex; see Figure~\ref{fig-locconv}(a). A~less trivial example is illustrated
in Figure~\ref{fig-locconv}(b). Given $c \in (0,\infty )$, the sequence
of Erd\H{o}s--R\'{e}nyi graphs $G(n,c/n)$ converges in distribution  in the local weak sense to the Galton--Watson
(GW) tree with offspring distribution given by the Poisson($c$) distribution. The latter is an example
of a unimodular Galton--Watson (UGW) tree, which is defined as follows. Given a
probability distribution $\rho $ on ${\mathbb{N}}\cup \{0\}$ that has finite
nonzero first moment, that is, satisfies
$0 < \sum_{k \in {\mathbb{N}}} k \rho (k) < \infty$, the random tree UGW($\rho$) has a root whose neighborhood
size is distributed according to~$\rho$.  The neighbors of the vertices are referred to as the offspring of
the root and form the  first generation of the tree.  Recursively, for  $n \geq 1$, each vertex in the nth
generation of the tree has an independent random number of offspring (equivalently, neighbors that
are further away from the root than itself) with  distribution $\widehat{\rho }$ 
\begin{equation}
\label{def2-hatrho} \hat{\rho }(k) = \frac{(k+1)\rho (k+1)}{\sum_{n \in {\mathbb{N}}}n\rho (n)}, \quad k \in {
\mathbb{N}}\cup \{0\}.
\end{equation}
The $(n+1)$th generation of the tree is comprised of all offspring of vertices in the $n$th generation.
It is easy to verify that if $\rho$ is a Poisson distribution, then
$\hat{\rho } = \rho $. Hence, a Galton--Watson tree with a Poisson($c$)
offspring distribution is in fact a UGW (Poisson($c$)) distribution.  Another special case is the $\kappa $-regular tree,
for $\kappa \geq 2$, which is given by
$\mathbb{T}_{\kappa }:=$ UGW ($\delta _{\kappa }$).
UGW trees are in a sense canonical
objects since they arise as local weak limits of many sparse
random graph sequences including Erd\H{o}s--R\'{e}nyi graphs,
configuration models and preferential attachment graphs;
see \cite[Section 2.2.4]{LacRamWuLWC21} for further discussion
 of these examples.\looseness=-1
\end{remark}

To extend this notion of convergence to graphs that are not necessarily connected, given an (unrooted) graph $G = (V,E)$ and a vertex $v \in V$, define $\mathsf{C}_{v}(G) \in \mathcal{G}_{*}$ to be the isomorphism
class of the connected component of $G$ that contains $v$, with $v$ as
its root. Furthermore, when $G$ is finite, we let $o$ denote a
random vertex of $G$ chosen uniformly from the set $V$, in which case $\mathsf{C}_{o}(G)$ denotes the connected
component of that random vertex.

\begin{definition}
\label{def-lwc}
A sequence of finite (random) graphs $\{G_{n}\}$ is said to
\emph{converge in distribution in the local weak sense} to $G$ if
\begin{align}
\lim_{n\to \infty }{\mathbb{E}} \biggl[\frac{1}{ \llvert G_{n} \rrvert }\sum
_{v \in G_{n}} f\bigl(\mathsf{C}_{v}(G_{n})
\bigr) \biggr] = {\mathbb{E}}\bigl[f(G)\bigr], \quad \forall f \in
C_{b}(\mathcal{G}_{*}). \label{def:localweakdist}
\end{align}
A sequence of finite (random) graphs $\{G_{n}\}$ is said to
\emph{converge in probability in the local weak sense} to $G$ if for every
$\varepsilon > 0$,
\begin{align}
\lim_{n\to \infty } \mathbb{P} \biggl( \biggl\llvert
\frac{1}{ \llvert G_{n} \rrvert }\sum_{v
\in G_{n}} f\bigl(
\mathsf{C}_{v}(G_{n})\bigr) - {\mathbb{E}}\bigl[f(G)
\bigr] \biggr\rrvert > \varepsilon \biggr) \rightarrow 0, \quad \forall f \in
C_{b}( \mathcal{G}_{*}). \label{def:localweakprob}
\end{align}
\end{definition}

Analogously, given a marked (unroooted, not necessarily connected) graph
$(G,x)$, $\mathsf{C}_{v}((G,x))$ 
denotes the connected component of $G$ containing $v$, 
with $v$ as its root and with the corresponding marks.
The notions of convergence
in distribution and in probability in the local weak sense for marked graphs
are defined in an exactly analogous fashion as Definition~\ref{def-lwc}, with $\mathsf{C}_{v}((G_n,x^n))$,  $\mathsf{C}_{v}((G_,x))$ and $(G,x)$ in place of
$\mathsf{C}_{v}(G_n)$,  $\mathsf{C}_{v}(G)$ and $G$, respectively.
For both unmarked and marked graphs,  convergence
in probability clearly implies convergence in distribution. We will use the same
notation for graphs and their isomorphism classes and often omit the root
from the notation and simply refer to $G \in \mathcal{G}_{*}$ rather than
$(G,o) \in \mathcal{G}_{*}$.

\begin{remark}
Given any  sequence
of random graphs $\{G_{n}\}_{n \in {\mathbb{N}}}$
that converges (either in distribution or in probability) in the local weak sense to a limit graph $G$, if
$x^{n} = (x^{n}_{v})_{v \in G_{n}}$ are i.i.d. marks on some Polish space $\mathcal{S}$ with the same distribution
irrespective of $n$, then it is easy to show that the marked graph sequence
$\{(G_{n},x^{n})\}_{n \in {\mathbb{N}}}$ also converges (in the same local
weak sense as the unmarked counterparts) to $(G,x)$, where
$x = (x_{v})_{v \in G}$ is i.i.d. with the same distribution. In fact,
as shown in \cite[Proposition 2.16]{LacRamWuLWC21}, convergence  of the marked graph sequence holds for the larger class of possibly dependent marks distributed according to a  Gibbs measure on the graph with respect to a fixed pairwise interaction functional.
\end{remark}

\subsection{Hydrodynamic limits} 
\label{subs-hydro}

\subsubsection{Interacting Diffusion processes}
\label{188.sec3.2.1}

We now address Q1--Q3 raised at the beginning of Section~\ref{sec-hydro}. The first result below states that if a sequence of graphs
marked with initial conditions converges (either in probability or in distribution)
in the local weak sense to a limit graph, then the graphs marked with the
trajectories that solve the corresponding SDE also converge to the limit graph in the same sense. The result also
characterizes the limit of the global empirical measure under suitable
conditions.

\begin{theorem}
\label{th-hydrodiff}
Suppose Assumption~\ref{ass-MFdiff} holds, and
the sequence
$\{(G_{n},\xi ^{n})\}_{n \in {\mathbb{N}}}$ of (not necessarily connected,
finite) random marked graphs converges in distribution in the local weak
sense to a $\mathcal{G}_{*}[B_{r}(\mathbb{R})]$-valued limit
$(G,\xi )$ for some $r > 0$.  Also,
for each $n \in {\mathbb{N}}$, let
$X^{G_{n},\xi ^{n}}$ be the solution to the SDE~\eqref{SDE-diff} with initial  data $(G_{n},\xi ^{n})$ and let $X^{G,\xi }$ be the unique weak solution to the SDE~\eqref{SDE-diff} on the limit graph $(G,\xi )$.
Then $\{(G_{n},X^{G_{n},\xi ^{n}})\}_{n \in {\mathbb{N}}}$ converges in distribution
in the local weak sense to the $\mathcal{G}_{*}[\mathcal{C}]$-valued element
$(G,X^{G,\xi })$.
In particular,
$\{X^{G_{n},\xi ^{n}}_{o_{n}}\}_{n \in {\mathbb{N}}}$ converges weakly to
$X^{G,\xi }_{o}$.  Moreover, if
$\{(G_{n},\xi ^{n})\}_{n \in {\mathbb{N}}}$ converges in probability in
the local weak sense to $(G,\xi )$, then
$\{(G_{n},X^{G_{n},\xi ^{n}})\}_{n \in {\mathbb{N}}}$ also converges in
probability in the local weak sense to $(G,X^{G,\xi })$ and additionally,
$\{\mu ^{G_{n},\xi ^{n}}\}_{n \in {\mathbb{N}}}$ converges weakly to the
law of $X_{o}^{G,\xi }$.
\end{theorem}

This result follows from \cite[Theorems 3.3 and 3.7]{LacRamWuLWC21},
which establish this result for more general, possibly non-Markovian diffusive
dynamics. A version of the first assertion of the above theorem was also
established for a slightly different class of interacting diffusions in
\cite{OliReiSto20}. Theorem~\ref{th-hydrodiff} can be seen as establishing
continuity in the local weak topology of the dynamics with respect to the initial
data, comprising the graph marked with initial conditions. It also provides
conditions under which the empirical measure can be shown to have a deterministic limit (equivalently, hydrodynamic limit) that additionally coincides with the limit law of the root particle, thus answering in the affirmative Q2 and Q3 at the beginning of Section~\ref{sec-hydro}.  As discussed earlier, on complete graphs the analogous phenomena holds due to asymptotic independence of the trajectories of any two particles. In contrast, in the sparse regime, neighboring particles remain dependent in the limit.  Instead, the proof relies on showing that the trajectories on finite neighborhoods
of two independent vertices, both chosen uniformly at random from the graph become asymptotically independent in the limit. The latter property  relies on a certain correlation decay property of the dynamics in the spirit of~\eqref{cordecay}; see \cite[Lemma 5.2]{LacRamWuLWC21} for details.

However, it should be emphasized that in the sparse regime the deterministic hydrodynamic limit
result holds \emph{only} when the initial data converges in the stronger sense of  convergence \emph{in probability} in the local weak sense.
Indeed, as shown in
\cite[Theorems 3.9 and 6.4]{LacRamWuLWC21}, the limiting empirical measure
can be stochastic when the initial data only converges in distribution
in the local weak sense. In particular, fix $c \in (0,\infty )$ and suppose
$\tilde{G}_{n}$ is the graph obtained by taking the \emph{connected component}
of a vertex chosen uniformly at random from the Erd\H{o}s--R\'{e}nyi graph
$G(n,c/n)$, and setting the root to be that chosen vertex. Also, let the
initial conditions $\xi ^{n} = \{\xi ^{n}_{v}\}_{v \in G_{n}}$ be i.i.d.
with common distribution $\gamma $, and let $\tilde{\xi }^{n}$ denote the
restriction of the initial conditions to $\bar{G}_{n}$. Then both
$\{(G_{n},\xi ^{n})\}_{n \in {\mathbb{N}}}$ and
$\{(\tilde{G}_{n}, \tilde{\xi }^{n})\}_{n \in {\mathbb{N}}}$ converge in
distribution in the local weak sense to $(G,\xi )$, where
$\xi = (\xi _{v})_{v \in V}$ is i.i.d. with distribution $\gamma $, and
$G = \mathcal{T}:= \mathrm{UGW}(\mathrm{Poiss}(c))$. However, whereas
$\mu ^{G_{n},\xi }$ converges weakly to
$\mathcal{L}(X_{o}^{\mathcal{T},\xi })$, the law of the dynamics at the
root vertex in $\mathcal{T}$,
$\mu ^{\tilde{G}_{n}, \tilde{\xi }^{n}}$ converges weakly to the following
\emph{random limit} $ \tilde{\mu }^{\mathcal{T},\xi }$ given by
\begin{equation*}
\tilde{\mu }^{\mathcal{T},\xi } := %\left \{
\begin{cases}
\mu ^{\mathcal{T},\xi } & \mbox{on the event }  \llvert \mathcal{T} \rrvert  < \infty ,
\\
\mathrm{Law} \bigl(X_{o}^{\mathcal{T},\xi } \mid  \llvert \mathcal{T} \rrvert  = \infty \bigr) &
\mbox{on the event }  \llvert \mathcal{T} \rrvert  = \infty .
\end{cases} 
\end{equation*}
This limit is truly stochastic  because, as is well
known from the elementary theory of branching processes
\cite{AthreyaNeybook}, there is always a positive probability for the UGW tree $\mathcal{T}$ to be finite (and there is a positive probability of $\mathcal{T}$  being infinite only when $c > 1$).  Furthermore, there also exist examples that show that even when
the limiting empirical measure is deterministic, it need not coincide with
the law of the root particle in the limit. For example, this can occur
if the graph itself is not homogeneous and the root is not chosen uniformly
at random from the vertices of the graph (see, e.g.,
\cite[Section 3.6]{LacRamWuLWC21}). The above discussion shows that the existence and
nature of the hydrodynamic limit is far more subtle in the sparse regime
than in the case of complete or dense graphs, even for diffusive dynamics.

\subsubsection{Interacting jump processes} 
\label{subs-hydrojump_}

In the setting of jump processes, additional subtleties arise. Whereas in the diffusion
setting, Assumption~\ref{ass-MFdiff} ensures that the drift of a particle
at a node $v$ experiences only an $O(1/|\partial v|)$ effect when there
is a perturbation in the state of a neighboring particle, an analogous assumption would be too stringent to cover
most jump models of interest on sparse graphs. In the latter case,
the effect of a neighboring particle on the jump intensity at a vertex
either remains constant or $O(1)$ as in the voter model, or  grows with the degree of the vertex in many other models, including the contact process (see \cite{Lig99}).  It is precisely to accommodate such a dependence that the jump intensity $\bar{r}_{j}$ in~\eqref{SDE-jump} is expressed as a function of the \emph{unnormalized} sum
of the Dirac masses at neighboring states introduced in~\eqref{theta_v}, rather than the normalized empirical measure. For hydrodynamic limits on sparse graphs, it will suffice to impose the following  mild assumption on the rates.

\begin{assumption}
\label{ass-MFjump}
For every $T > 0$, suppose there exist constants
$C_{k,T}, k \in {\mathbb{N}}$, such that $k \mapsto C_{k,T}$ is nondecreasing
and for every $j \in \mathcal{J}$,
$\sup_{x \in \mathcal{X}, t \in [0,T]} \bar{r}_{j}(t,x,\theta )
\leq C_{\theta (\mathcal{X}),T}$.
\end{assumption}

Note that in the jump SDE \eqref{SDE-jump} describing the dynamics, the third argument $\theta$ of $\bar{r}_{j}$  is equal to the unnormalized empirical measure  of the states of the neighbors of a vertex, as defined in \eqref{theta_v}. Thus,  $\theta(\mathcal{X})$ irepresents the degree of the vertex and
 Assumption \ref{ass-MFjump} allows the uniform bound on the jump rates at a vertex to grow with the degree of a vertex.

We now state an analog of Theorem~\ref{th-hydrodiff} for jump diffusions.
Recall the definition of a UGW tree given in Remark~\ref{rem-locconv}.

\begin{theorem}
\label{th-hydrojump}
Suppose Assumption~\ref{ass-MFjump} holds, and the sequence
$\{(G_{n},\xi ^{n})\}_{n \in {\mathbb{N}}}$ of (not necessarily connected,
finite) random rooted marked graphs converges in probability in the local weak
sense to a limit $(G,\xi )$, where $G$ is a UGW($\rho $) tree with
$\rho $ having finite, strictly positive first and second moments. For
each $n \in {\mathbb{N}}$, let $X^{G_{n},\xi ^{n}}$ be the solution to the
jump SDE~\eqref{SDE-jump} with initial data $(G_{n},\xi ^{n})$, and let $X^{G,\xi }$ be the unique strong solution to the
SDE~\eqref{SDE-jump} on the limit graph $(G,\xi )$. Then
$\{(G_{n},X^{G_{n},\xi ^{n}})\}_{n \in {\mathbb{N}}}$ converges in probability
in the local weak sense to the $\mathcal{G}_{*}[\mathcal{D}]$-valued element
$(G,X^{G,\xi })$.  Furthermore,
$\{\mu ^{G_{n},\xi ^{n}}\}_{n \in {\mathbb{N}}}$ converges weakly to the
law of $X_{o}^{G,\xi }$.
\end{theorem}

This theorem follows from more general results established in
\cite[Theorem 4.8 and Corollary 5.16]{GanRam-Hydro22}.  As  in the diffusion case, the proof of
the theorem involves establishing continuity properties of the dynamics with respect to the graph and initial condition as well as a correlation decay property.  However, the proofs of these properties are considerably more involved than in the diffusion case due to the weaker conditions imposed on the jump intensities in Assumption~\ref{ass-MFjump}. For one,  in the jump
setting even well-posedness of the particle system on an infinite random
graph of unbounded degree is not automatic.  As shown in
\cite[Appendix B]{GanRam-Hydro22}, there exist examples of simple jump
particle systems with uniformly bounded jump rate functions that can have
multiple solutions on certain graphs with exponential growth. To quote
Liggett \cite{Lig72}, ``Given an intuitive description of the behavior
of the particles, it is often not clear whether or not there exists a ...
process which corresponds to that description. Therefore it is important
to find conditions under which infinite particle systems exist.'' On finite
graphs, there are only finitely number of jumps in any bounded interval, and the process remains constant between jumps.
Thus, one can simply order the jumps and define the process recursively.
The problem in the infinite graph setting is that one cannot always identify
a ``first'' jump. In \cite{Lig72}, Liggett used an analytical construction
invoking the theory of semi-groups to establish a general existence theorem
for the law of Markovian particle systems on infinite graphs with quite general (not necessarily finite-range) interactions.  In the context of nearest-neighbor interactions on lattices,
an alternative, probabilistic approach was used to establish well-posedness of Markovian interacting particle systems in Harris \cite{Har72,Har74}.
However,  both approaches seem to only apply to graphs with finite maximal degree.
 On the other hand, graphs of particular
interest like the UGW(\textrm{Poisson}($c$)) tree discussed above (see Remark~\ref{rem-locconv}) have unbounded degrees.  Under
Assumption~\ref{ass-MFjump},  well-posedness of (possibly non-Markovian) interacting jump processes was established in \cite[Theorem 4.3]{GanRam-Hydro22} for
a large class of possibly random graphs that satisfy a certain ``finite
dissociability'' property almost surely, and this property
was shown to hold for UGW trees in
\cite[Corollary 5.16]{GanRam-Hydro22}. The proof of Theorem~\ref{th-hydrojump} then
follows on combining this well-posedness result with continuity properties of the dynamics (with respect to the initial data) and a correlation decay property (established in
\cite[Proposition 6.8]{GanRam-Hydro22} and
\cite[Theorem~4.9]{GanRam-Hydro22}, respectively).

\section{Marginal dynamics on trees} 
\label{sec-marg}

The hydrodynamic limit result reduces the characterization of the limit law of $X^{G_{n},\xi ^{n}}_{o_{n}}$ to the understanding of the marginal law $X^{G,\xi }_{o}$ of the root dynamics on the \emph{infinite} limit graph $G$. Given that local weak limits of many random graphs are trees (see Remark~\ref{rem-locconv}),  we focus here on understanding marginal dynamics on (random) trees. The evolution of the root $X^{G_{n},\xi ^{n}}_{o_{n}}$ in~\eqref{SDE-diff2} or~\eqref{SDE-jump} is driven by the \emph{local} neighborhood empirical measure $\mu ^{G_{n},\xi ^{n}}_{o_{n}}$ or its unnormalized counterpart
$\theta ^{G_{n},\xi ^{n}}_{o_{n}}$. In the complete (or sufficiently dense)
graph case, in the limit as the number of particles goes to infinity, the local empirical measure coincides with the
global empirical measure. Thus, in this case the hydrodynamic limit  yields
an autonomous characterization of the limit marginal dynamics. In contrast,
when the graph sequence $G_{n}$ is sparse,  neighboring vertices remain strongly correlated, the local neighborhood
empirical measure remains stochastic and thus the hydrodynamic limit results
in Section~\ref{sec-hydro} are not adequate to provide an autonomous characterization
of the marginal dynamics.

Instead, we adopt a different perspective, which is better suited to the analysis of
large collections of dependent random elements. As mentioned in Section~\ref{subs-back}, as a first step we try to identify the conditional independence
structures in such random variables. To this end, we identify a certain
Markov random field (MRF) property for the trajectories of
$X^{G,\xi } = \{X^{G,\xi }_{v}, v \in G\}$ in Section~\ref{subs-MRF} below.
We then describe how to exploit this property, along with filtering results
from stochastic analysis and symmetry properties of the graph, to identify
an autonomously defined ``local equation'' satisfied by marginal dynamics on $\mathrm{cl}_{o}$, the \emph{root and its neighborhood}.
We do this first for diffusions on the line in Section~\ref{subs-derivation}, then for diffusions on UGW trees in Section~\ref{subs-margUGW}, and finally for jump processes in Section~\ref{subs-margjump}.  Unlike in the mean-field case, consideration of the marginal at the root \emph{and} its neighborhood (rather than just at the root), appears necessary in order to obtain an autonomous characterization.  This is also necessary in order to capture correlations between neighboring vertices, which do not vanish in the sparse regime.
Further discussion of the local equation  is given in Sections~\ref{subs-derivation}-\ref{subs-margjump}. But it is worth noting here that since the graphs we consider are locally finite, the neighorhood of the root is (almost surely) finite. Thus,
the  local equation describes the evolution of an (almost surely) finite number of  interacting  particles.  When combined with the convergence results of Theorems \ref{th-hydrodiff} and \ref{th-hydrojump},  this finite-dimensional interacting process serves as an approximation for the marginals of (possibly non-Markovian) interacting processes with an arbitrarily large number of particles, and thus constitutes a significant dimension reduction.

\subsection{A Markov random field property} 
\label{subs-MRF}

We first introduce the definition of an MRF.

\begin{definition}
\label{def-MRFs}
Fix a measurable space $\mathcal{Y}$, and a (possibly infinite, but locally
finite) graph $G = (V,E)$. A random element $Y = (Y_{v})_{v \in V}$ is said to be a \emph{(first-order) MRF} (abbreviated as MRF) on
$\mathcal{Y}^{V}$ with respect to $G$ if for every \emph{finite} set
$A \subset V$, $Y_{A}$ is conditionally independent of
$Y_{(A \cup \partial A)^{c}}$ given $Y_{\partial A}$, which we denote as
\begin{equation}
\label{def-1MRF} Y_{A} \indep Y_{(A \cup \partial A)^{c}} |
Y_{\partial A}.
\end{equation}
On the other hand, $Y$  is said to be a \emph{(first-order)   global MRF}  if~\eqref{def-1MRF} holds for
all $A \subset V$, possibly infinite.
Furthermore $Y$  is said to be a \emph{(first-order)
  semi-global MRF} (abbreviated as SGMRF) if~\eqref{def-1MRF} holds for all $A \subset V$ such that
$\partial A$ is finite.
Furthermore, $Y = (Y_{v})_{v \in V}$ is said to be a \emph{second-order MRF} or 2-MRF (respectively,
second-order SGMRF or 2-SGMRF) with respect to $G$ if it is an MRF (respectively,
SGMRF) with respect to the square graph $G^{2} = (V,E^{2})$, where
$E^{2}$ contains $E$ as well as vertex pairs that are a distance two apart
in $G$.
In all cases, we will say
$\nu \in \mathcal{P}(\mathcal{Y}^{V})$
exhibits a certain MRF property whenever some
$\mathcal{Y}^{V})$-valued
random element $Y$ with law $\nu$ satisfies that MRF property.
\end{definition}

The SGMRF property, introduced in \cite{GanRam-Hydro22}, can be viewed as a generalization of tree-indexed Markov
chains \cite[Chapter 12]{Geo11} to general graphs, and is clearly strictly stronger
than the MRF property. For any $t > 0$, and collections of paths
$x = (x_{v})_{v \in V}$ (either in $\mathcal{C}$ or $\mathcal{D}$), let

\begin{equation}
\label{not-paths} x_{v}[t] := \bigl(x_{v}(s)
\bigr)_{s \in [0,t]} \quad \mbox{and} \quad x_{v}[t) :=
\bigl(x_{v}(s)\bigr)_{s \in [0,t)}
\end{equation} 
represent the trajectory of $x_{v}$ in the intervals $[0,t]$ and
$[0,t)$, respectively, and for any subset $A \subset V$, let
$x_{A}[t] := (x_{v}[t])_{v \in A}$ and
$x_{A}[t) := (x_{v}[t))_{v \in A}$.   Certain MRF
properties are preserved under the evolution of  interacting processes, in a sense made precise in the following theorem.

\begin{theorem}
\label{th-MRF}
Let $G = (V,E)$ be a (deterministic) graph with uniformly bounded  degree or the almost sure realization of a UGW tree.  Suppose the $\mathbb{R}^V$-valued element
$\xi = (\xi _{v})_{v \in V}$ forms a $2$-MRF (or $2$-SGMRF) with respect to $G$,  Assumption~\ref{ass-MFdiff} holds and $X^{G,\xi }$ is the unique solution to the diffusive SDE~\eqref{SDE-diff}. Then the $\mathcal{C}$-valued trajectories
$X^{G,\xi } = (X^{G,\xi }_{v})_{v \in V}$ also form a $2$-MRF  (respectively, $2$-SGMRF) with respect to $G$.   On the other hand, suppose  Assumption~\ref{ass-MFjump} holds,
$\xi = (\xi _{v})_{v \in V}$  is a $\mathcal{X}^V$-valued random element that forms a $2$-MRF (or $2$-SGMRF) with respect to $G$,  and
$X^{G,\xi}$ is the unique solution to the jump
SDE~\eqref{SDE-jump}.
Then $X^{G,\xi} = (X^{G,\xi }_{v})_{v \in V}$ also forms a $2$-MRF (respectively,
$2$-MRF) on $G$ in $\mathcal{D}$. In both cases, the same assertions also hold with
$X^{G,\xi }$ replaced with $X^{G,\xi }[t]$ or $X^{G,\xi }[t)$ for any $t \geq 0$.
\end{theorem}

The discussion in Section~\ref{subs-derivation} provides insight into why only the second-order, and not in general the first-order, MRF property
is preserved by the dynamics.
The preservation of the $2$-MRF property for diffusions for graphs with bounded degree follows from \cite[Theorem 2.7]{LacRamWuMRF21}.  The proof proceeds by first establishing the result on finite graphs by appealing to Girsanov's theorem and the Gibbs--Markov theorem \cite[Theorem 2.30]{Geo11} (also often referred to as the Hammersley--Clifford
theorem), and then suitably approximating infinite systems by a sequence of finite systems.
The proof of preservation of both the 2-MRF and 2-SGMRF properties for jump processes in \cite[Theorem 3.7]{GanRam-MRF22}
follows a rather different approach. It exploits a certain
duality between marginals of the interacting system and nonexplosive point
processes to directly establish an infinite-dimensional Girsanov theorem,
obviating the need for any approximation arguments. This approach also
allows more general initial conditions that can incorporate infinite histories,
which is required to characterize solutions to the local equation described in Section \ref{subs-margjump} as flows on a suitable path space \cite{GanRam-DET22}.
The result for diffusions in \cite{LacRamWuMRF21} can be
generalized in a similar fashion using the approach developed
in~\cite{GanRam-MRF22}.

Prior work on such questions has
largely focused on interacting diffusions,  specifically  characterizing them as Gibbs
measures on path space in order to construct weak solutions to infinite-dimensional
SDEs.   Deuschel \cite{Deu87} initiated this perspective for diffusions with drifts of gradient type. Although not explicitly stated, the 2MRF property is    implicit
in his proof of existence of the weak solution, which relies on estimates of Dobrushin's contraction
coefficient that crucially require additional smoothness and boundedness
properties of the drift. Cattiaux, Roelly, and Zessin
\cite{CatRoeZes96} relaxed the boundedness condition to allow Markovian,
Malliavin differential drifts, using a variational characterization and
an integration-by-parts formula. Subsequent works
\cite{MinRoeZes00,PraRoe04} used a cluster expansion method that applies
to systems obtained as small perturbations of non-interacting systems.  Dereudre and Roelly \cite{DereudreRoelly17} established
Gibbsian properties of paths of  interacting one-dimensional diffusions on $\mathbb{Z}^{m}$ with (possibly history-dependent) drift having sublinear
growth using specific entropy, but this crucially requires shift-invariant initial conditions.
In another direction, several other works have considered the MRF (or Gibbsian)
nature of marginals rather than of paths, both in the diffusion and jump process
contexts \cite{RedRoeRus10,RoeRus14,Vanetal02,Kul19,KisKul20}, but preservation
of this property holds in general only for sufficiently small time horizons
or interaction strengths.  Furthermore,
none of this work seems to have considered the SGMRF property, which is crucial for the derivation of the local equation, as elaborated in Sections \ref{subs-derivation}-\ref{subs-margjump} below.

\subsection{Outline of derivation of the local equation for diffusions on the line} 
\label{subs-derivation}

We now describe how the $2$-SGMRF property of Theorem~\ref{th-MRF} can be used
to obtain an autonomous characterization of the marginal dynamics of the
root neighborhood on the $2$-regular tree $\mathbb{T}_{2}$. For simplicity, we identify  $\mathbb{T}_{2}$ with $\mathbb{Z}$ and identify the root
$o$ with $0$.  Additionally,  rather than consider the general form in~\eqref{SDE-diff2} with $G = \mathbb{Z}$, we focus on the special case of pairwise interacting diffusions in~\eqref{SDE-diff}, but without  time dependence in the drift:
\begin{align}
dX_{v}(t) = \frac{1}{2} \bigl[ \beta
\bigl(X_{v}(t),X_{v+1}(t)\bigr) + \beta
\bigl(X_{v}(t), X_{v-1}(t)\bigr) \bigr] dt +
dW_{v}(t), \quad v \in \mathbb{Z}, \label{def:T2-particlesystem}
\end{align}
where we have dropped the superscripts denoting graph dependence for notational conciseness.
We also assume that $(X_{v}(0))_{v \in \mathbb{Z}}$ is a shift-invariant 2-SGMRF.

Given the above setup,
our goal is to understand the marginal dynamics
$X_{\{-1,0,1\}} = (X_{-1}, X_{0}, X_{1})$ of the root \emph{and its neighborhood}.
The characterization of this marginal via the local equation entails four key ingredients, which we elaborate upon below.

\vskip.3cm

\noindent
(i) \emph{A mimicking theorem.}
By~\eqref{def:T2-particlesystem},
we can rewrite the dynamics of the marginal $X_{\{-1,0,1\}}$ of interest as follows:
\begin{equation}
\label{SDE-reduced} dX_{v} (t) = b_{v}(t,X) dt +
dW_{v}(t), \quad v \in \{-1,0,1\},
\end{equation}
where for $v \in \{-1,0,1\}$ and $X = (X_{v})_{v \in \mathbb{Z}}$,
\begin{equation}
\label{SDE-drift} b_{v}(t,X) := \frac{1}{2} \bigl[\beta
\bigl(X_{v}(t),X_{v+1}(t)\bigr) + \beta
\bigl(X_{v}(t), X_{v-1}(t)\bigr) \bigr], \quad v \in
\{-1,0,1\}.
\end{equation}
Let $(\Omega , \mathcal{F}, \mathbb{F}, \mathbb{P})$ denote the filtered
space that supports the $\mathbb{F}$-adapted process $X$.  For the
root node, observe that at time $t$, the drift $b_{0}$ of $X_0$ depends on $X$ only through
$X_{\{-1,0,1\}}(t)$. However, at time $t$ the drift $b_{1}$ of $X_1$ depends on $X_{2}(t)$ and likewise the drift $b_{-1}$ of node $-1$ depends on $X_{-2}(t)$. Since $2$ and $-2$ do not lie in the closure $\{-1,0,1\}$ of the root, the system of equations  \eqref{SDE-reduced} is not autonomous since its drift
at time $t$ depends on random elements beyond $X_{\{-1,0,1\}}(t)$. Nevertheless,~\eqref{SDE-reduced} and~\eqref{SDE-drift} together show that
$X_{\{-1,0,1\}}$ is what is known as an It\^{o} process, which means that
its drift $(b_{-1}, b_{0}, b_{1})$ is $\mathbb{F}$-progressively measurable
(as a consequence of~\eqref{SDE-drift}, the continuity of $\beta $ and
the fact that $X$ is an $\mathbb{F}$-adapted continuous process). Therefore,
one can appeal to a ``mimicking'' theorem for It\^{o} processes from filtering
theory (see \cite{liptser-shiryaev} or
\cite[Appendix A]{LacRamWuLE21}), which allows one
to express  $X_{\{-1,0,1\}}$ as the solution to an SDE
 whose drift at time $t$ is a functional only of the past of $X_{\{-1,0,1\}}$ up to time $t$, rather than an arbitrary
 $\mathbb{F}$-adapted process.  Then the mimicking theorem allows one to conclude that
(by extending the probability space if necessary) there exist independent
Brownian motions
$(\widetilde{W}_{-1},\widetilde{W}_{0},\widetilde{W}_{1})$ on the extended probability space such that
$X=(X_{-1},X_{0},X_{1})$ satisfies
\begin{equation}
\label{xv-tildeb}
X_{v} (t) = \tilde{b}_{v} (t, X) dt + d
\widetilde{W}_{v} (t), \quad v \in \{-1,0,1\},
\end{equation}
where for $v \in \{-1,0,1\}$,
$\tilde{b}_{v}:[0,\infty ) \times \mathcal{C}^{\{-1,0,1\}} \mapsto
\mathbb{R}^{d}$ is a progressively measurable version of the conditional
expectation:
\begin{equation}
\label{def-tildebv} \tilde{b}_{v} (t, x) := \frac{1}{2} {
\mathbb{E}} \bigl[ \beta \bigl( X_{v}(t), X_{v+1}(t)
\bigr) + \beta \bigl(X_{v}(t),X_{v-1}(t)\bigr) \, \big| \,
X_{\{-1,0,1\}}[t] = x[t] \bigr].
\end{equation}
Recall from~\eqref{not-paths} that
$x[t] = (x(s))_{s\in [0,t]}$. Clearly, the conditioning does not alter
the drift coefficient for the root particle, which remains the same as
in the original system~\eqref{def:T2-particlesystem}:
\begin{equation}
\label{tilde-bo} \tilde{b}_{0} (t,x) = \frac{1}{2} \bigl[
\beta \bigl(x_{0}(t), x_{1}(t)\bigr) + \beta
\bigl(x_{0}(t),x_{-1}(t)\bigr) \bigr].
\end{equation}
However, $\tilde{b}_{1}$ and $\tilde{b}_{-1}$ will be altered by the conditioning.  We now see how the MRF property can be used to  simplify the expression for $\tilde{b}_{1}$ and $\tilde{b}_{-1}$.

\vskip.3cm

\noindent
(ii) \emph{Markov random field structure.} To compute the drifts
$\tilde{b}_{1}$ and $\tilde{b}_{-1}$ and provide a self-contained description
of the law of the dynamics of $X_{\{-1,0,1\}}$, one needs to be able to
express the conditional law of $X_{2}(t)$ and $X_{-2}(t)$ given the past
$X_{\{-1,0,1\}}[t]$ in terms of the (joint) law of $X_{\{-1,0,1\}}$ or
preferably, in terms of the law of $X_{\{-1,0,1\}}[t]$ to get a nonanticipative
description of the dynamics. To this end, we invoke the property from Theorem~\ref{th-MRF} that the trajectories
$(X_{i}[t])_{i \in \mathbb{Z}}$ up to time $t$ form a 2-SGMRF or second-order
Markov chain in $\mathbb{Z}$:
\begin{align}
\bigl(X_{j}[t]\bigr)_{j < i} \indep \bigl(X_{j}[t]
\bigr)_{j > i + 1} \, \big|\, \bigl(X_{i}[t],X_{i+1}[t]
\bigr), \ \ \forall i \in \mathbb{Z}. \label{def:T2-condind}
\end{align}

Before we use this property, let us consider the corresponding
first-order property: namely to ask whether we should in fact expect that for every $t > 0$,
\begin{equation*}
\bigl(X_{j}[t]\bigr)_{j < i} \indep \bigl(X_{j}[t]
\bigr)_{j > i} \,\big|\, X_{i}[t], \quad \forall i \in
\mathbb{Z}.
\end{equation*}
One may attempt to bolster this hypothesis by reasoning that conditioned
on $X_{i}[t] = \psi $, $X_{i-1}[t]$ and $X_{i+1}[t]$ become decoupled and
satisfy the following SDE: for $s \in [0,t]$,
\begin{align*}
dX_{i-1}(s) & = \frac{1}{2}\bigl[\beta
\bigl(X_{i-1}(s), \psi (s)\bigr) + \beta \bigl(X_{i-1}(s),
X_{i-2}(s)\bigr)\bigr] ds + dW_{i-1}(s),
\\
dX_{i+1} (s) & = \frac{1}{2}\bigl[\beta
\bigl(X_{i+1}(s), X_{i+2}(s)\bigr) + \beta
\bigl(X_{i+1}(t), \psi (s)\bigr)\bigr] ds + dW_{i+1}(s),
\end{align*}
where $W_{i-1}$ and $W_{i+1}$ are independent Brownian motions. However,
a more careful inspection would reveal that such a reasoning is fallacious
because the evolution of $X_{i}$, and thus the random element
$X_{i}[t]$, directly depends on the states
$(X_{i-1}(s), X_{i+1}(s))_{s \in [0,t]}$, which are in turn dependent on
the Brownian motions $W_{i-1}$ and $W_{i+1}$. Thus, conditioning on $X_{i}[t] = \psi $ causes the driving Brownian motions $W_{i-1}[t]$ and $W_{i+1}[t]$ to become correlated. Thus, under this conditioning, $X_{i-1}$ and $X_{i+1}$ are not independent and do not follow the above SDE driven by independent Brownian motions on $[0,t]$.  However,~\eqref{def:T2-condind} shows that by conditioning on both $X_{i}[t]$ and
$X_{i+1}[t]$, the driving noise processes $W_{i-1}$ and $W_{i+2}$ remain decoupled (i.e., independent), although the conditioning does alter the distributions of
$W_{i-1}$ and $W_{i+2}$; they are no longer Brownian motions or even martingales.

Returning to the simplification of the expression for the drifts
$\tilde{b}_1$ and $\tilde{b}_2$ given in~\eqref{def-tildebv}, note that since the relation in~\eqref{def:T2-condind} implies that
$X_{2}(t)$ is independent of $X_{-1}[t]$ when conditioned on
$X_{\{0,1\}}[t]$, the drift $\tilde{b}_{1}$ in~\eqref{def-tildebv} can be rewritten as
\begin{equation}
\begin{split} \tilde{b}_{1} (t, x) & :=
\frac{1}{2} {\mathbb{E}} \bigl[ \beta \bigl( X_{1}(t),
X_{0}(t)\bigr) + \beta \bigl(X_{1}(t),X_{2}(t)
\bigr) \, \big| \, X_{
\{-1,0,1\}}[t] = x[t] \bigr]
\\
& = \frac{1}{2} \beta \bigl( x_{1}(t),
x_{0}(t)\bigr) + \frac{1}{2} {\mathbb{E}} \bigl[\beta
\bigl(X_{1}(t),X_{2}(t)\bigr) \, \big| \,
X_{\{0,1\}}[t] = x_{\{0,1
\}}[t] \bigr]. \label{tilde-b-1}
\end{split} %
\end{equation}
By the same reasoning, an analogous expression holds for $\tilde{b}_{-1}$.

\vskip.3cm

\noindent
(iii) \emph{Symmetry considerations.} Despite the simplification of the last section, the second term on the right-hand
side of~\eqref{tilde-b-1} still involves $X_{2}$, and thus has  not been written purely in terms of  $X_{-1,0,1}[t]$ and its law.  However, it can
be rewritten in this form by exploiting the shift-variance of the particle
system on $\mathbb{Z}$ (since this is true of both the initial condition
and the dynamics). More precisely,  the fact that
$(X_{0}, X_{1}, X_{2})$ has the same distribution as
$(X_{-1}, X_{0}, X_{1})$ allows us to conclude that
\begin{equation*}
\tilde{b}_{1} (t, x) = \frac{1}{2} \beta \bigl(
x_{1}(t), x_{0}(t)\bigr) + \frac{1}{2} {
\mathbb{E}} \bigl[\beta \bigl(X_{0}(t),X_{1}(t)\bigr) \,
\big| \, X_{
\{-1,0\}}[t] = x_{\{0,1\}}[t] \bigr].
\end{equation*}
Along with the analogous expression for $\tilde{b}_{-1}$,
and equations~\eqref{xv-tildeb},~\eqref{def-tildebv}, and ~\eqref{tilde-bo},  this shows that
\begin{equation}
\begin{aligned} dX_{0}(t) &= \frac{1}{2} \bigl[
\beta \bigl(X_{0}(t), X_{1}(t)\bigr)dt + \beta
\bigl(X_{0}(t), X_{-1}(t)\bigr) \bigr] dt + d
\widetilde{W}_{0}(t),
\\
dX_{k}(t) &= \frac{1}{2} \bigl[ \beta
\bigl(X_{k}(t), X_{0}(t)\bigr) + \tilde{\gamma
}_{t}(X_{k},X_{0}) \bigr] dt + d
\widetilde{W}_{k}(t), \quad k \in \{-1,1\}, {.}
\end{aligned}
\label{def:intro:T2localeq}
\end{equation}
where $\widetilde{W}_{-1}, \widetilde{W}_{0}$ and
$\widetilde{W}_{1}$ are independent $d$-dimensional Brownian motions, and
\begin{align}
\tilde{\gamma }_{t}(x,y) :={\mathbb{E}} \bigl[ \beta
\bigl(X_{0}(t),X_{1}(t)\bigr) \, \big| \, (X_{0},X_{-1})[t]
= (x,y)[t] \bigr], \quad (x,y) \in \mathcal{C}^{2}. \label{def:intro:gammatil}
\end{align}
Modulo some additional technical (measurability and integrability) conditions,
this identifies the form of the \emph{local equation} satisfied by the
$\{-1,0,1\}$ marginal dynamics on the line (see
\cite[Definition 3.5 with $\kappa = 2$]{LacRamWuLE21} for a complete definition).

Observe that even though the original system~\eqref{def:T2-particlesystem} describes a (linear) Markov process, its
marginal $X_{\{-1,0,1\}}$, characterized by the local equation system~\eqref{def:intro:T2localeq}--~\eqref{def:intro:gammatil}, is a \emph{nonlinear}, describes a nonlinear
\emph{non-Markovian} process since the drift functional $\tilde{\gamma }_{t}$ depends on both the history of $X_{\{-1,0,1\}}[t]$  up to time $t$ and its law.  However, the structure of
$\tilde{\gamma }_{t}$ ensures that the coupled system~\eqref{def:intro:T2localeq} no longer depends on values of the process  $X_v$ for $v \notin \{-1,0,1\}$, and is thus autonomously defined.

\vskip.3cm
\noindent
(iv) \emph{Uniqueness of solutions to the local equation.} The above argument shows that the law of the marginal solves the local equation system~\eqref{def:intro:T2localeq} and~\eqref{def:intro:gammatil}.  To complete the autonomous \emph{characterization}
of the marginal it only remains to show that the law of the marginal of $X_{-1,0,1}$ is the \emph{unique}
solution to the local equation system (or rather, its complete specification as stated in
\cite[Definition 3.5 with $\kappa = 2$]{LacRamWuLE21}). The methods described
in Section~\ref{subs-MFdiff} to prove well-posedness of nonlinear Markov
processes (which characterize the limit marginal law of a node in the complete
graph case) all run into difficulties here due to the path-dependence and,
more importantly, the nonlinearity occurring through dependence on conditional
laws, which are less regular. Nevertheless, it is possible to prove well-posedness
using other approaches, entailing relative entropy estimates and symmetry
properties, or via a correspondence with the infinite particle system; further
details can be found in \cite[Sections 4.3.1 and 4.2]{LacRamWuLE21}.

The local equation on the root neighborhood of a $\kappa$-regular tree $\mathbb{T}_{\kappa }$, with $\kappa \geq 3$, can be derived in a   manner similar to the case of $\mathbb{T}_{2}$, once again invoking the mimicking theorem and Theorem~\ref{th-MRF},
but now exploiting the additional ``rotational and transational'' symmetries
arising from the automorphism groups of $\mathbb{T}_{\kappa }$, in place of just the translational and reflection symmetries of
$\mathbb{T}_2.$  However, the full expression of the local equation is omitted as it is a special
case of the UGW tree discussed in the next section, whose analysis is more subtle.

\subsection{Local equations for diffusions on unimodular Galton--Watson trees} 
\label{subs-margUGW}

Let $\rho $ be a probability distribution on
${\mathbb{N}}\cup \{0\}$ satisfying
$\sum_{k\in {\mathbb{N}}} k \rho _{k} < \infty $, and let
$\mathcal{T}$ be a UGW($\rho $) tree as in Remark~\ref{rem-locconv}. Again,
we would like to describe the marginal dynamics of the particle process
on the root node and its (random) neighborhood.
As elucidated in the last section, the main ingredients in the derivation
of the local equation on the line are a mimicking theorem, a conditional
independence property, symmetry considerations and, finally, well-posedness of the local equation.
The mimicking theorem can be applied without change also on
$\mathcal{T}$. However, the MRF property in Theorem~\ref{th-MRF} applies to deterministic graphs and is thus not sufficient.  Instead, one needs an \emph{annealed version}, that is, one that also takes into account the random structure of the tree $\mathcal{T}$.  For any $t > 0$, one has to show that (on the event the root is not isolated)
for any child $k$ of the root, conditioned on the trajectories
$(X^{\mathcal{T}}_{o}[t], X_{\kappa }[t])$, the trajectories
$X_{\mathcal{T}_{k}}^{\mathcal{T}}[t]$ of particles on the subtree
$\mathcal{T}_{k}$ rooted at $k$ are independent of the trajectories
$X_{\mathrm{cl}_{\mathcal{T}}(o)}[t]$ on the root and its neighborhood,
and, moreover, that the conditional law of
$X_{\mathcal{T}_{k}}^{\mathcal{T}}[t]$ given
$X^{\mathcal{T}}_{\{o,k\}}[t]$ does not depend on $k$ (see
\cite[Proposition 3.17]{LacRamWuLE21}). In the case when
$\rho = \delta _{\kappa }$, this would follow from Theorem~\ref{th-MRF} and homogeneity of the dynamics, but in the general case,
one also has to account for the randomness of $\mathcal{T}_{k}$ and the
root neighborhood.

Furthermore, the symmetry properties are now considerably more subtle -- the appropriate notion here being \emph{unimodularity}. Unimodularity can be viewed as the analog on an infinite graph of the property on the finite graph  that the root is uniformly distributed on the graph.
Since the latter statement is not well defined on an infinite graph, this is instead phrased in terms of a ``mass transport principle'' that on finite
graphs is equivalent to having a uniformly distributed root.  The definition of  unimodularity involves the space $\mathcal{G}_{**}$ of isomorphism
classes of doubly rooted graphs, which is defined as follows, in a fashion
analogous to $\mathcal{G}_{*}$. A \emph{doubly rooted graph}
$(G,o,o')$ is a rooted graph $(G,o)$ with an additional distinguished vertex
$o'$ (which may equal $o$). Two doubly rooted graphs
$(G_{i},o_{i},o_{i}')$ are isomorphic if there is an isomorphism from
$(G_{1},o_{1})$ to $(G_{2},o_{2})$ which also maps $o_{1}'$ to
$o_{2}'$. We write $\mathcal{G}_{**}$ for the set of isomorphism classes
of doubly rooted graphs. A double rooted marked graph is defined in the
obvious way, and $\mathcal{G}_{**}[\mathcal{S}]$ denotes the set of isomorphism
classes of doubly rooted marked graphs. The space
$\mathcal{G}_{**}[\mathcal{S}]$ is equipped with its Borel $\sigma $-algebra.

\begin{definition}
\label{def-unimodular}
For a metric space $\mathcal{S}$, a $\mathcal{G}_{*}[\mathcal{S}]$-valued
random variable $(G,o,S)$ is said to be \emph{unimodular} if the following
\emph{mass-transport principle} holds: for every (nonnegative) bounded Borel
measurable function
$F : \mathcal{G}_{**}[\mathcal{S}] \rightarrow \mathbb{R}_{+}$, 
\begin{equation*}
\mathbb{E} \biggl[ \sum_{o' \in G} F
(G,S,o,o') \biggr] = \mathbb{E} \biggl[ \sum
_{o' \in G} F(G,S,o',o) \biggr].
\end{equation*} 
\end{definition}
Combining all these properties it was
 shown in \cite{LacRamWuLE21} that the marginal of the particle system $X^{\mathcal{T}}$ on the closure $\bar{\mathcal{T}}$ of the root neighborhood of the UGW ($\rho $) tree $\mathcal{T}$ can be characterized by a local equation that has a similar form to~\eqref{def:intro:T2localeq}, except
that $\tilde{\gamma }_{t}$ is now a reweighted version of the conditional
expectation of the drift that takes into account the structure of the tree:
on the event that the root is not isolated,
\begin{equation*}
\gamma _{t}(X_{o},X_{1}) = 2
\frac{{\mathbb{E}} [\alpha _{t} \beta (X_{o},X_{\partial _{\bar{\mathbb{T}}}(o)}) \mid
X_{o}[t], \, X_{1}[t] ]}{{\mathbb{E}} [ \alpha _{t}\mid X_{o}[t], \, X_{1}[t] ]},
\end{equation*}
where
$\alpha _{t} := |\partial _{\bar{\mathcal{T}}}(o)|/(1+ \widehat{C}_{1})$,
with $\hat{C}_{1}$ being a random variable distributed according to~$\hat{\rho }$
(representing the number of offspring of a child of the root)
that is independent of the root neighborhood structure, the initial conditions
and the driving Brownian motions, and the factor of $2$ in the expression
arises just to compensate for the $1/2$ that arises in~\eqref{def:intro:T2localeq}. This extra weighting by $\alpha _{t}$ arises due to the unimodularity property and significantly complicates the proof of well-posedness of the local equation.

\subsection{Marginal dynamics for pure jump processes} 
\label{subs-margjump}

As in the last section, let $\mathcal{T} = (V,E)$ be the UGW($\rho$) tree,  and given initial conditions
 $\xi = (\xi)_{v \in V}$, let $X = X^{{\mathcal T}, \xi}$  be  the solution of the jump SDE~\eqref{SDE-jump} with $G = \mathcal{T}$, and also denote $\theta_v = \theta_{v}^{G,\xi }$.   We are once again interested in obtaining an autonomous characterization of the marginal  law of $X$ on
   the root and its neighborhood in terms of a corresponding  local equation.    The derivation of this local equation follows the same broad  outline that was used in the case of diffusions, although the justification of each step require substantially different arguments.  For simplicity, we flesh out a few details in the special case when
$\mathcal{T}$ is the rooted $\kappa $-regular tree $\mathbb{T}_{\kappa }$ for $\kappa \geq 2$.  Let $\tilde{\mathbb{T}}_\kappa$ denote the subtree of $\mathbb{T}_\kappa$ consisting of the root and its neighborhood.  Then,
  on verifying certain technical conditions, one can first appeal to filtering results  for point processes (see, e.g.,  \cite{Daley-VereJones}) to establish an analogous  mimicking theorem for jump processes.  Specifically, we use the latter result to argue that $X$ restricted to $\tilde{\mathbb{T}}_\kappa$ can be expressed (on a possibly extended probability space) as the solution to the following functional jump SDE:  for $t \geq 0$,
   \begin{equation}
\label{SDE-jump-mimicking} X_{v}(t) = X_{v}(0) +
\sum_{j \in \mathcal{J}} j \int_{(0,t]
\times \mathbb{R}_{+}}
\mathbb{I}_{\{r\leq \tilde{r}_{j}(s,X)\}} \tilde{N}_{v}(ds,dr),
\quad v \in \tilde{\mathbb{T}}_\kappa,
\end{equation}
where $(\tilde{N}_{v})_{v \in V}$, are i.i.d. Poisson random measures on $\mathbb{R}_{+}^{2}$ with intensity measure $\Leb ^{2}$,  and $\tilde{r}_j: \mathbb{R}_+ \times \mathcal{D}^{\kappa+1} \to \mathbb{R}_+$ is a predictable version of the conditional expectation
  \[   \tilde{r}_j^v (t, x) = \mathbb{E}\left[ \bar{r}_j^v(t, X_v(t-),\theta_v(t-))| X_{\tilde{\mathbb{T}}_\kappa}[t) = x[t)\right],
  \]
  for a.s. every $x$.

  Combining this with  the 2-SGMRF property for $X$ established in Theorem \ref{th-MRF} for interacting jump processes with  $G = \mathbb{T}_\kappa$,  and invoking the symmetries of the law of $X$ with respect to the automorphisms of $\mathbb{T}_\kappa$, the equation~\eqref{SDE-jump-mimicking} can  be further simplified into an autonomous local equation of the following form (see \cite{GanRam-DET22}):
\begin{align}
X_{o}(t) &= X_{o}(0) + \sum
_{j\in \mathcal{J}}\int_{[0,t]
\times \mathbb{R}_{+}}
\mathbb{I}_{\{r \leq \bar{r}_{j} (s, X_{o}(s-),
\theta _{o}(s-))\}} \tilde{N}_{0}(ds,dr),
\nonumber \\[-6pt]
\\[-6pt]
X_{k}(t) &= X_{k}(0) + \sum
_{j\in \mathcal{J}}\int_{[0,t]\times
\mathbb{R}_{+}} \mathbb{I}_{\{r \leq \tilde{\gamma}_{j} (s, X_{k}, X_{o})
\}}
\tilde{N}_{k}(ds,dr), \quad k = 1, \ldots , \kappa ,
\nonumber
\end{align} 
where recall $\theta _{o}(s-) = \sum_{k=1}^{\kappa }\delta _{X_{k}(s-)}$ and
$\tilde{\gamma}: \mathbb{R}_{+} \times \mathcal{D}^{2} \mapsto \mathbb{R}_{+}$
is defined by 
\begin{equation}
\label{eq-tilder} \tilde{\gamma}_{j} (t,x,y) := {\mathbb{E}} [
\bar{r}_{j} \bigl(t,X_{o}(t), \theta _{o}(t)
\bigr)\mid (X_{o},X_{1})[t) = (x,y)[t) ].
\end{equation}
Here, once again, we have omitted various measurability and other technical conditions
required to define a solution to the local equation, referring the reader to \cite{GanRam-DET22} for full details.  As in the case of diffusions, it is evident from the local equation
the marginal dynamics is non-Markovian even when  the original dynamics
in~\eqref{SDE-jump} is, and it is also nonlinear in the sense that the evolution of the process depends on its own law.  The local equations identify precisely the nature of this nonlinear non-Markovian dynamics.
As in the diffusion case, it is also possible to define analogous local equations describing marginal dynamics on the
UGW tree, which rely on a more involved (annealed) SGMRF property and a more complicated proof of well-posedness of the local equation.

\subsection{Generalizations and approximations} 
\label{subs-gen}

For both diffusive and jump dyncamis, one could also consider non-Markovian interacting
processes. For example, consider solutions to the SDE~\eqref{SDE-diff2} in which the drift $b$ at the vertex $v$ is replaced
by a suitably regular nonanticipative functional
$F_{v}: \mathbb{R}_{+} \times \mathcal{C}^{V} \to \mathbb{R}$.
For example, consider $F_{v}(t,x) = b(t,x_{v} (t-\tau ), \mu _{v}(t-\tau ))$, with
$\mu _{v}(s) = \frac{1}{\partial v} \sum_{u \sim v} \delta _{x_{u}(s)}$
for some $\tau > 0$.  Or likewise, consider solutions to the jump SDE~\eqref{SDE-jump} in which the jump rate $\bar{r}_{j}^v$ at  vertex $v$ is
replaced with a suitable predictable functional of the paths such as
$F_{j,v}(t,x) = \bar{r}_{j}(t,x_{v}(t-\tau ),\theta _{v}(t-\tau ))$, where
$\theta _{v}(s) = \sum_{u \sim v} \delta _{x_{u}(s)}$. It is not too difficult
to see from the discussions of the derivation of the local equation given in Sections~\ref{subs-derivation}-\ref{subs-margjump}  that even
in the non-Markovian setting, under suitable regularity conditions, the marginal law of the process on the root and its neighborhood
could  be characterized by an analogous local equation. Indeed,
the frameworks in both
\cite{LacRamWuLWC21,LacRamWuMRF21,LacRamWuLE21} and
\cite{GanRam-Hydro22,GanRam-MRF22,GanRam-DET22} allow for quite general
non-Markovian dynamics. Furthermore, one can also allow more general initial
conditions that are not necessarily i.i.d. but form a $2$-SGMRF
and (in the non-Markovian setting) incorporate histories of the process
up to time $0$. The latter is useful for studying flow properties
of the local equation dynamics and gaining insight into stationary measures for the local equations
\cite{GanRam-DET22}. Furthermore, the framework in
\cite{GanRam-Hydro22} can also be used to handle interacting particle systems
on directed graphs (see \cite[Remark 2.2]{GanRam-Hydro22}).

In the case when the full system is non-Markovian, the local equation yields
a significant dimension reduction parallel to that achieved in mean-field
limits, since it approximates the marginal of a non-Markovian system on
an arbitrarily large high-dimensional (random) graph by a nonlinear non-Markovian
process of a fixed  finite (average) dimension. On the other hand, when one seeks
to use the local equations to approximate the marginal of a Markovian system,
since the local equation is still non-Markovian. Thus, in terms of computing or simulating the process,  there is a tradeoff between the size of the time interval one is interested in and the size (or number of particles) in the original
Markovian system. It is thus natural to ask if any further principled approximations are
possible in that case to  make the local equations more analytically and
computationally tractable even over long time intervals.

\begin{figure}[b!] 
\includegraphics{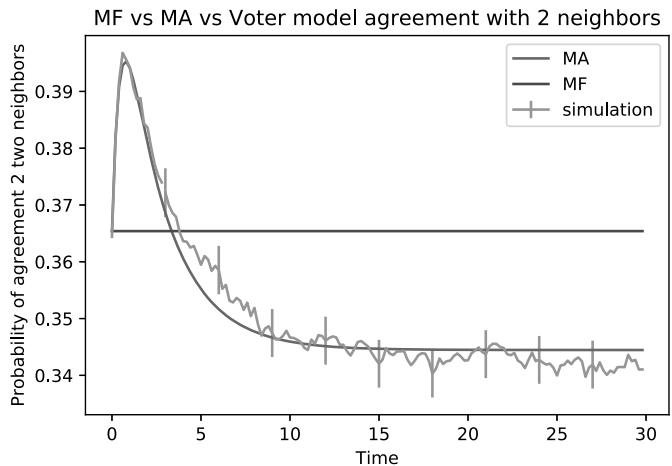}
\caption{Comparing simulations, the mean-field approximation and a Markovian version of the local equation}
\label{fig-votermodel2} 
\end{figure}

Recall from Figure~\ref{fig-votermodel} that for the voter model on the
$3$-regular tree $\mathbb{T}_{3}$ (truncated after $9$ generations) the mean-field
approximation for the probability of agreement of the root with precisely
two of its neighbors was rather inaccurate. Figure~\ref{fig-votermodel2} plots,
for the same tree and parameters, an \emph{ad hoc}  Markovianization of the local
equation, wherein $\tilde{\gamma}_{j}$ in \eqref{eq-tilder} is replaced with a modified state-dependent version $\bar{\gamma}_j: \mathbb{R}_+ \times \mathcal{X}^2 \to \mathbb{R}_+$ given by
\[
\bar{\gamma}_{j} (t,x,y) := {\mathbb{E}} [
\bar{r}_{j} \bigl(t,X_{o}(t), \theta _{o}(t)
\bigr)\mid (X_{o},X_{1})[t-) = (x,y)[t-) ].
\]
The good agreement of the simulation with the
Markovian version of the local equation in Figure~\ref{fig-votermodel2} suggests that such a Markovian local equation  may serve as a good approximation for several models. This motivates a more rigorous investigation of the accuracty of the Markovian local equation for various classes of models and derivation of rigorous error
bounds between the laws of the solution to the Markovianized local equation and the original local equation.

\section{Open questions} 
\label{sec-open}

This article describes the first reduced dimension characterization of
marginal dynamics on sparse random graphs, thereby resolving an open question raised
in \cite{DelGiaLuc16} in the context of interacting diffusions. A plethora of open questions remain, related to
the structure of the solution to the local equation such as ergodic properties, as well as
theoretical guarantees for developing more computationally tractable principled approximations, and also applications
(see, for example,  \cite{Ram-QUESTA22}). A few open questions are
listed below.

\subsection{Long-time behavior and invariant measures for the local equations}
% \end{description}
%
The results thus far have focused on transient dynamics over finite time
intervals. There are several open questions about equilibrium behavior
and long-time behavior.
\vspace*{-3pt}
\begin{enumerate}
\item[{Q1.}] Can one can establish general conditions for existence and
uniqueness of stationary or invariant measures for the local equation?

\item[{Q2.}] Can one identify which of these stationary measures correspond
to a stationary measure of the evolution on the corresponding infinite graph? Can the local equation
be used to identify phase transitions (i.e., identify parameters for existence of multiple stationary distributions on random trees)?

\item[{Q3.}] Can one study ergodic properties and use the local equation
to sample from marginal stationary distributions on random graphs?
\vspace*{-3pt}
\end{enumerate}
Recent work \cite{LacZha21,GanRam-DET22} has analyzed some stationarity
properties for interacting systems indexed by a regular tree. The work
\cite{LacZha21} studies limits of systems of diffusions with gradient drift
that have an explicit Gibbs measure (i.e., $1$-MRF) as the unique invariant
measure on any finite graph, and relates the stationary distribution of
the full system to that of a modified local equation. Continuous-state
MRFs on infinite trees can also be studied via recursions (see, e.g.,
\cite{RozBook12,GamRam19}). On the other hand, the work
\cite{GanRam-DET22} studies jump processes with possibly nonreversible
dynamics.  Also, in \cite{CocRam23} the long-time behavior of the Susceptible-Infected-Recovered
(SIR) interacting particle system on UGW trees is shown to be succinctly characterized in terms of
a fixed point equation. \\

% \vspace{-0.1in}

\subsection{Refined convergence results} 

%\begin{description}
%\item[
% \noindent 
% \textbf{{B. Refined convergence results}} 
%\end{description}

% \vspace{-0.0001in}

A natural question is whether one can obtain more refined convergence results,
that provide concentration results and rates of convergence, as well as
a characterization of fluctuation and large deviations from the hydrodynamic
limit. Large deviations principles also provide an alternative way of characterizing
hydrodynamic limits.

% \vspace{-2.8in}

\begin{enumerate}
\item[{Q4.}] Can one establish large deviation principles and concentration
results for interacting particle systems on sparse random graphs?
\vspace*{-3pt}
\end{enumerate}
%

% \vspace{-0.8in}

Such results have been obtained for weakly interacting particle systems
on complete and dense graphs (see, e.g.,
\cite{DawGar87,BudDupFis12,OliRei18,CopDieGia18,DelLacRam20,Ram21,BalPerRei22}).

 \vspace{-2.8in}

% \begin{description}

\subsection{Analytic characterizations}

% \noindent 
% \textbf{C. Analytic characterizations} 
 %  \vspace*{-3pt}
% \end{description}
%

% \vspace{-2.8in}

In the mean-field setting, the corresponding nonlinear process and its
stationary distribution can also be described by\vadjust{\vspace*{-6pt}\eject} nonlinear PDEs (in the
diffusive case) or nonlinear integrodifferential equations (in the jump
case).
\vspace*{-3pt}
\begin{enumerate}
\item[{Q5.}] Is it possible to develop a corresponding theory for these
new types of path-dependent nonlinear equations that involve conditional
laws? Also, can one determine when the marginal laws are absolutely continuous
with respect to Lebesgue measure?
\vspace*{-3pt}
\end{enumerate}

\subsection{From interacting particle systems to games}

Mean-field approximations have been used to study not only interacting
particle systems but also games where where strategic agents control their
dynamics to maximize an objective function. When the dynamics and objective
functions are symmetric, a limit problem called the mean-field game has
shown to provide tractable approximations to Nash equilibria in finite-agent
games, which are notoriously hard to compute (see \cite{Del-MFG21} for
surveys on different aspects of mean-field games).
\vspace*{-3pt}
\begin{enumerate} 
\item[{Q6.}] Can one establish limit theorems for Nash equilibria of games
with a large number of agents in which the interaction network of agents
is sparse rather than the complete graph? While there have been several
recent results looking at mean-field games on networks with nodes whose
degrees diverge to infinity, there are only a few works studying this on
graphs with uniformly bounded degree (see \cite{LacSor21} for the study
of linear--quadratic games and the works
\cite{DetFouIch18,IchFenFou21a,IchFenFou21b} for games on directed graphs).  
\vspace*{-3pt}
\end{enumerate}

\vspace{0.1in}

\noindent 
{\bf Acknowledgments. } 
  This work was partially supported by~the Vannevar Bush Faculty Fellowship
ONR-N0014-21-1-2887, Army Research Office grant 911NF2010133 and the National
Science Foundation grant DMS-1954351.  
I would also like to acknowledge several 
collaborators J. Cocomello,  A. Ganguly, D. Lacker, and R. Wu for various joint works on this subject, 
 thank A. Ganguly for generating Figures~\ref{fig-graphs}, \ref{fig-votermodel},
and \ref{fig-votermodel2}, and thank J. Cocomello and K. Hu for comments.

%

%%%% BackMatter %%%
%\begin{ack}
%\end{ack}

%\vfill\eject

%\begin{funding}
%\end{funding}

% \bibliographystyle{plain}
% \bibliography{ICM.bib}

% \end{document}

% structpyb loaded by romualda, 2022-03-10 09:35:56

%
\end{document}